\documentclass[11pt,twoside,a4paper]{article}
\usepackage{amsfonts}
\usepackage{mathrsfs}
\usepackage{mathrsfs}
\usepackage{mathrsfs}
\usepackage{mathrsfs}
\usepackage{mathrsfs}
\usepackage{mathrsfs}
\usepackage{mathrsfs}

\usepackage{amsfonts}
\usepackage{amsmath,amssymb}
\usepackage{mathrsfs}
\usepackage{hyperref}

\newtheorem{thm}{Theorem}[section]
\newtheorem{cor}[thm]{Corollary}
\newtheorem{lem}[thm]{Lemma}
\newtheorem{prop}[thm]{Proposition}

\numberwithin{equation}{section}\allowdisplaybreaks

\def\le{\leqslant}
\def\ge{\geqslant}
\def\leq{\leqslant}
\def\geq{\geqslant}

\def\Real{{\mathbb{R}}}

\begin{document}

\title{ {\bf \Large Sharp global well-posedness for non-elliptic derivative Schr\"odinger equations with small rough data}}

\author{\textsc{\normalsize Baoxiang Wang\footnote{
 Email: wbx@math.pku.edu.cn}}       \\
 $^{*}$\textit{\footnotesize {LMAM, School of Mathematical Sciences, Peking University and BICMR,
Beijing 100871, China}}  }
\date{}
\maketitle

\begin{minipage}{13.5cm}
\footnotesize \bf Abstract. \rm We show the sharp global well posedness for the Cauchy problem for the cubic (quartic) non-elliptic derivative Schr\"odinger equations with  small rough data in modulation spaces $M^s_{2,1}(\mathbb{R}^n)$ for $n\ge 3$ ($n= 2$).  In 2D cubic case,
 using the Gabor frame, we get some time-global dispersive estimates for the Schr\"odinger semi-group in anisotropic Lebesgue spaces, which include  a time-global maximal function estimate in the space $L^2_{x_1}L^\infty_{x_2,t}$.  By resorting to the smooth effect estimate together with the dispersive estimates in anisotropic Lebesgue spaces, we show that the cubic hyperbolic derivative NLS in 2D has a unique global solution if the initial data in Feichtinger-Segal algebra or in weighted Sobolev spaces are sufficiently small.

\vspace{10pt}



\end{minipage}

\section{Introduction}

In this paper we consider the Cauchy problem for the derivative nonlinear Schr\"odinger equation (DNLS):
\begin{align}
{\rm i} u_t -  \Delta_\pm u = F(u, \bar{u}, \nabla u, \nabla\bar{u}), \quad u(0,x)=u_0(x),\label{dnls}
\end{align}
where $\Delta_\pm = \varepsilon_1\partial^2_{x_1} +...+ \varepsilon_n \partial^2_{x_n}$, $\varepsilon_i\in \{1,\ -1\}$ for $i=1,...,n$,  $\nabla u= (u_{x_1},..., u_{x_n})$,
\begin{align}
F(z)= F (z_{1},...,
z_{2n+2})=\sum_{m+1 \leq\mid  \beta \mid
<\infty}c_{\beta}z^{\beta}, \ \ m\ge 2,  \quad c_{\beta}\in \mathbb{C}
 \label{nonlpoly}
\end{align}
and $|c_\beta| \le C^{|\beta|}$ for $\beta= (\beta_1,...,\beta_{2n+2})$.  A special case of \eqref{dnls} is the following
\begin{align}
{\rm i} u_t -  \Delta_\pm u = \overrightarrow{\lambda} \cdot\nabla(|u|^{2\kappa} u)
  + \mu |u|^{2\nu}u, \quad u(0,x)=u_0(x),
\label{nonl}
\end{align}
$\vec{\lambda}  \in \mathbb{C}^n$, $\mu\in \mathbb{C}$ and $\kappa, \nu \in \mathbb{N}$. It is known that  \eqref{nonl} ($\kappa=1$) arises from the strongly interacting many-body
systems near criticality as recently described in terms of nonlinear
dynamics \cite{ClTu90, DiTu89, TuDi89}, where anisotropic interactions are manifested by the presence of
the non-elliptic case, as well as additional residual terms which
involve cross derivatives of the independent variables.    Another typical example is the Schr\"odinger map equation
\begin{align}
{\rm i} u_t -  \Delta_\pm u     = \frac{2\bar{u}}{1+|u|^2} (\varepsilon_1 u^2_{x_1}+...+ \varepsilon_n u^2_{x_n}),  \quad u(0,x)=u_0(x), \label{Schmap}
 \end{align}
which is an equivalent form of the non-elliptic Schr\"odinger map
\begin{align}
 s_t   =s \times \Delta_\pm s,    \
 \  s( 0,x)=s_0(x).   \label{Hei}
\end{align}
where $s=(s_1,s_2,s_3)$,  $s: \mathbb{R} \times \mathbb{R}^n \to \mathbb{S}^2$ is a real valued  map  of $(t, x_1, ..., x_n)
$.   Indeed, if $u$ satisfies  \eqref{Schmap}, taking
$$
s=\left(\frac{2Re \ u}{1+|u|^2}, \  \frac{2Im \ u}{1+|u|^2}, \ \frac{1-|u|^2}{1+|u|^2}\right),
$$
we see that $s$ is the solution of \eqref{Hei}. Conversely, taking $u$ as the stereographic projection of $s$ defined by
$$
u= \frac{s_1+is_2}{1+s_3},
$$
we see that \eqref{Hei} reduces to \eqref{Schmap}. A large amount of work has been devoted to the study of the elliptic Schr\"odinger map initial value problem ($\Delta_\pm=\Delta$) together with their generalizations \cite{BeIoKeTa08, DiWa01, MeRaRo11, So92, WaHaHu09, ZhGuTa91}.

If there is no condition on initial data, it is easy to give a finite-time blow up solution  of the derivative NLS and in the non-elliptic case, all of the blow up points can constitute a curve,  see Appendix C.

For the general equation \eqref{dnls} in the elliptic case with nonlinearity \eqref{nonlpoly},  the local and global well posedness were studied in \cite{Ch96, Kl82, KlPo83, OzZh08, Sh82}.  When the nonlinear term $F$ satisfies
an energy structure condition $Re \partial F /\partial (\nabla u) =0 $ and the initial data are sufficiently smooth in weighted Sobolev spaces,  Klainerman \cite{Kl82}, Shatah \cite{Sh82} and Klainerman and Ponce \cite{KlPo83} obtained the global existence of \eqref{dnls} in all spatial dimensions.   Chihara \cite{Ch96} considered the initial data in sufficiently smooth weighted Sobolev spaces and removed the energy structure condition $Re \partial F /\partial (\nabla u) =0 $ for $n\ge 3$ and only assume that cubic terms $F_1(z)$ in $F(z)$ is modulation homogeneous (i.e., $F_1(e^{{\rm i}\theta}z)= e^{{\rm i}\theta} F_1(z)$)  for $n=2$.  Ozawa and Zhai \cite{OzZh08} was able to consider the initial data in $H^s$ with $s>n/2+2$, $n\ge 3$ and $Re \partial F /\partial (\nabla u) =\nabla (\theta(|u|^2))$ for some real valued function $\theta\in C^2$  with $\theta(0) = 0$.

In the non-elliptic case, the smooth effect estimates seem to be useful tools for the well posedness of \eqref{dnls} with nonlinearity \eqref{nonlpoly}. Roughly speaking, this method relies upon the dispersive structure for the Schr\"odinger semi-group and the energy structure conditions for the nonlinear terms are not necessary for the local well posedness and for the global well posedness with small data.   By setting up the local smooth effects for the solutions of the
linear Schr\"odinger equation, Kenig, Ponce and Vega
\cite{KePoVe93, KePoVe98} were able to deal with the non-elliptical
case and they established the local well posedness of Eq.
\eqref{dnls} in $H^s$ with $s\gg n/2$. The global existence and scattering of solutions of \eqref{dnls} and \eqref{nonlpoly} with small data in modulation spaces $M^{5/2}_{2,1}$   (so in $H^s$ with $s>n/2+5/2$, $n\ge 3$) were recently obtained in \cite{WaHaHu09}.  Moreover, the results in \cite{WaHaHu09} contains non-elliptic Schr\"odinger map equation as a special case if $n\ge 3$.

In this paper we study the global well posedness of solutions of \eqref{nonl}, and \eqref{dnls} and \eqref{nonlpoly}, we show that \eqref{nonl}
in modulation spaces $M^s_{2,1}$ has a critical index $s=1/2\kappa$ for which \eqref{nonl} is globally well posed for small data if $s\ge 1/2\kappa$ and ill posed if $s<1/2\kappa$, $\kappa\ge 1$ for $n\ge 3$, $\kappa\ge 2$ for $n=2$ respectively.  Similarly, \eqref{dnls} with nonlinearity \eqref{nonlpoly} in modulation spaces $M^s_{2,1}$ has a critical index $s=1+1/m$ for which it is globally well posed with small data if $s\ge 1+1/m$ and ill posed  if $s<1+1/m$,  $m\ge 2$ for $n\ge 3$, $m\ge 3$ for $n=2$ respectively. In 2D case with $m=2$, and 1D case with $m=3$,
we show that \eqref{dnls} with nonlinearity \eqref{nonlpoly} is globally well posed for the small Cauchy data in
Feichtinger-Segal's algebra $M^{2}_{1,1}$.  On the basis of the Gabor frame expression for the initial data,  we will establish a class of linear estimates in anisotropic Lebesgue spaces,  these estimates together with the smooth effect estimates for the linear Schr\"odinger
equation  (cf. \cite{CoSa88, KePoVe93, KePoVe98, LiPo93, Sj87, WaHaHu09, Ve88}) and frequency-uniform
decomposition techniques yield the existence and uniqueness of global solutions for small initial data.

\subsection{Notation }
In the sequel $C$, $C_i$ will denote   universal positive constants which
can be different at each appearance. $x\lesssim y $  for $x$, $y
>0$  means that $x\le Cy$, and $x\sim y$ stands for $x\lesssim y $
and $y\lesssim x$; $x\vee y=\max(x,y)$. For any $p\in [1,\infty]$, $p'$ denotes the
dual number of $p$, i.e., $1/p+1/p'=1$.

  Let  $\mathscr{S} $ be
Schwartz space and $\mathscr{S}' $ be its dual space. All of the function spaces used in this paper are subspaces of $\mathscr{S}' $.  We will use the Lebesgue spaces $L^p:=
L^p (\mathbb{R}^n)$ with the norm $\|\cdot\|_p:=
\|\cdot\|_{L^p }$, the function spaces $L^q_{t}
L^p_x$ and $L^p_x L^q_{t}$ for which the norms are defined by:
\begin{align}
 \|f\|_{L^q_{t} L^p_x}=\left\|\|f\|_{L^p_x (\mathbb{R}^n)}\right\|_{L^q_t (\mathbb{R})}, \quad  \|f\|_{L^p_x L^q_{t }}=\left\|\|f\|_{L^q_t (\mathbb{R})}\right\|_{L^p_x(\mathbb{R}^n)}. \nonumber
\end{align}
We denote by $L_{x_i}^{p_1}L_{ (x_j)_{j\neq i},t
}^{p_2} $ for $(x_j)_{j\neq i} =(x_1,...,x_{i-1}, x_{i+1},...,x_n)$ and $L_{x_1}^{p_1}L_{ \bar{x},t
}^{p_2} $ for $\bar{x}=(x_2,...,x_n)$ the anisotropic Lebesgue
space for which the following norm is finite:
\begin{align}
\|f\|_{ L_{x_i}^{p_1}L_{ (x_j)_{j\neq i},t
}^{p_2} }  =\Big\|\|f\|_{L_{(x_j)_{j\neq i},  t}^{p_2}(\mathbb{R}^{n}) }\Big\|_{L_{x_i}^{p_1}(\mathbb{R}) } ,  \ \
 \|f\|_{L_{x_1}^{p_1}L_{\bar{x}, t}^{p_2} }   =\Big\|\|f\|_{L_{\bar{x}, t}^{p_2}(\mathbb{R}^{n}) }\Big\|_{L_{x_1}^{p_1}(\mathbb{R}) }.
\end{align}
Let $\mathscr{F}$ ($\mathscr{F}^{-1}$) be the (inverse) Fourier transform.  We will write
 $D_{x_i}^s = (-\partial_{x_i}^2)^{s/2}=
\mathscr{F}_{\xi_i}^{-1}|\xi_i|^s\mathscr{F}_{x_i}$ denotes the
partial Rieze potential in the $x_i$ direction.  The
homogeneous Sobolev space $\dot{H}^{s} $ is defined by
$(-\Delta)^{-s/2}L^2 $, $H^s= L^2\cap \dot H^s$.
Recall that the weighted Sobolev space $H^{s,b}(\mathbb{R}^n)$ is defined by
$$
\|u\|_{H^{s,b}}= \|\langle x\rangle^b \mathscr{F}^{-1}\langle \xi\rangle^s \mathscr{F} u\|_2.
$$

  Let   $ \{\sigma_k\}_{k\in \mathbb{Z}^n}$ be a smooth cut-off function  sequence satisfying
\begin{align}
  {\rm supp} \sigma_0  \subset [-1,1]^n, \ \ \sigma_k = \sigma(\cdot - k), \ \
    \sum_{k\in\mathbb{Z}^n}\sigma_k(\xi)\equiv 1. \label{denote}
\end{align}
We know that if $\{\eta_k\}_{k\in \mathbb{Z}}$ satisfies \eqref{denote} for $n=1$, then we see that
\begin{align}
   \sigma_k(\xi) = \eta_{k_1}(\xi_1)... \eta_{k_n}(\xi_n), \ \ k=(k_1,...,k_n)  \label{denote2}
\end{align}
satisfies \eqref{denote}.  We can define the frequency-uniform decomposition
operators $\Box_k$ as:
\begin{align}
\Box_k := \mathscr{F}^{-1}\sigma_k \mathscr{F}, \quad
k\in\mathbb{Z}^n,
\end{align}
and we write
\begin{align}
&\|f\|_{M_{p,q}^s(\mathbb{R}^n)}= \left(\sum_{k\in\mathbb{Z}^n} \langle k\rangle^{sq}
\|\Box_k f\|^q_{L^p (\mathbb{R}^n)} \right)^{1/q},
\end{align}
which is said to be a modulation space. Modulation spaces $M^s_{p,q}$ were introduced by Feichtinger  with the following equivalent norm (cf. \cite{Gr01, Gr06})
\begin{align} \label{moddef}
\|f\|_{M^{s}_{p,q}}=\left(\int_{\Real^{n}}\left(\int_{\Real^{n}}|V_{g}f(x,\omega)|^{p}dx\right)^{q/p}\langle
\omega\rangle^{sq}d\omega\right)^{1/q},
\end{align}
where, for a given window function $g\in \mathscr{S}$, $V_g$ is the short-time Fourier transform:
\begin{align*}
V_{g}f(x,\omega)=\int_{\Real^{n}}e^{-{\rm i} t
\omega}\overline{g(t-x)}f(t)dt.
\end{align*}
The behavior of the Schr\"odinger semi-group in modulation spaces are rather different from those in Lebesgue spaces. Indeed, it was shown in \cite{BeGrOkRo06, WaZhGu06, WaHe07} that
\begin{align*}
& \|e^{{\rm i}t \Delta_\pm} u_0\|_{M^s_{p,q}} \lesssim  (1+|t|)^{|n(1/2-1/p)|}\|  u_0\|_{M^s_{p,q}}, \ \ 1\le p,q \le \infty, \ s\in \mathbb{R}, \\
& \|e^{{\rm i}t \Delta_\pm} u_0\|_{M^s_{p,q}} \lesssim  (1+|t|)^{-n(1/2-1/p) }\|  u_0\|_{M^s_{p',q}}, \ p> 2, \ 1/p+1/p'=1, q\ge 1, s\in \mathbb{R},
\end{align*}
which means that Schr\"odinger semi-group is bounded in any $M^s_{p,q}$ and satisfies a truncated decay from $M^s_{p',q}$ to $M^s_{p,q}$.  It is known that  Schr\"odinger semi-group in Lebesgue spaces $L^{p}$ is bounded if and only if $p=2$ and the truncated decay from $L^{p'}$ to $L^p$ does not holds.

Considering the inclusions between modulation and Sobolev spaces, we have (cf. \cite{SuTo07, WaHu07})
\begin{align*}
H^{n/2+ s_+} \subset M^s_{2,1} \subset H^s \ \ {\rm (sharp \ inclusions)}.
\end{align*}
From scaling point of view, we see that the critical Sobolev space of \eqref{nonl} in the case $\mu=0, n\ge 2$ is $H^{n/2-1/2\kappa}$.  If we can solve \eqref{nonl} in $M^s_{2,1}$ with $s < n/2-1/2\kappa$, which contains a class of data in $L^2\setminus H^{n/2-1/2\kappa}$, then there is a class of $H^{n/2-1/2\kappa}$ supercritical data such that \eqref{nonl} is well posed.

The space $M^s_{1,1}$ with $s\ge 0$, so called Feichtinger algebra (or Feichtinger-Segal algebra), is one of the most important modulation spaces which enjoys the following interesting property (\cite{Fe81, FeWe06}): $\mathscr{F} \, (\mathscr{F}^{-1}):  M^0_{1,1}(\mathbb{R}^n) \to M^0_{1,1}(\mathbb{R}^n)$ is an isometric mapping.
On the other hand, it is known that $B^{s+n}_{1,1}(\mathbb{R}^n) \subset  M^s_{1,1}(\mathbb{R}^n) \subset  B^s_{1,1}(\mathbb{R}^n)$ are sharp inclusions (cf. \cite{SuTo07, To04, WaHu07})\footnote{$B^s_{p,q}$ denotes Besov space.}. $M^0_{1,1}$ can be regarded as an analogue  of Schwartz space which preserves Fourier transform but has little smoothness.

\subsection{Main results}

For convenience, we write for $n\ge 2$,
\begin{align}
  \|u\|^{\rm sm} = & \sum^n_{i=1} \sum_{k\in \mathbb{Z}^n, \, |k_i|\ge 20 \vee \max_{j\neq i}|k_j|} \langle k_i\rangle^{1/2+1/m} \|\Box_k u\|_{L^\infty_{x_i}L^2_{(x_j)_{j\neq i},t}} , \\
 \|u\|^{\rm max} = &  \sum^n_{i=1} \sum_{k\in \mathbb{Z}^n}  \|\Box_k u\|_{L^m_{x_i}L^\infty_{(x_j)_{j\neq i},t}} , \\
  \|u\|^{\rm str} =&  \sum_{k\in \mathbb{Z}^2} \langle k \rangle^{1/m}   \|\Box_k u\|_{L^\infty_{t}L^2_{x} \,  \cap \, L^{m+2}_{x,t}}.
\end{align}
We denote $\|u\|^{\Upsilon\cap \Gamma }=   \|u\|^{\Upsilon }+\|u\|^{ \Gamma}$.
\begin{thm}\label{th1deri}
Let $n\ge 2$, $\kappa, \nu \in \mathbb{N}$, $\kappa\ge 2$ for $n=2$. Suppose that  $u_0\in M^{1/2\kappa}_{2,1}$ and there exists a suitably small $\delta>0$ such that $\|u_0\|_{M^{1/2\kappa}_{2,1}} \le \delta$. Then  \eqref{nonl}
has a unique solution $u \in C(\mathbb{R}, M^{1/2\kappa}_{2,1})  \cap X^{1/2\kappa}$, where
 \begin{align}
X^{1/2\kappa} = \left\{u\in \mathscr{S}':  \| u\|^{\rm sm \, \cap \, max  \cap \, str} \lesssim \delta \right\}, \ \ m=2\kappa. \label{X1d}
\end{align}
Moreover, if $s<1/2\kappa$, then \eqref{nonl} with $\mu =0$ is ill-posed in the sense that the solution map $u_0\to u$ is not $C^{2\kappa+1}$ from $M^s_{2,1}$ to $C([0,T]; M^s_{2,1})$ for any $T>0$.
\end{thm}

\begin{thm}\label{thmderig}
Let $n\ge 3$, $m\ge 2$.  Suppose that $u_0\in M^{1+1/m}_{2,1}$ and there exists a suitably small $\delta>0$ such that $\|u_0\|_{M^{1+1/m}_{2,1}} \le \delta$. Then \eqref{dnls} has a unique solution $u \in C(\mathbb{R}, M^{1+1/m}_{2,1})  \cap X^{1+1/m}$, where
 \begin{align}
X^{1+1/m}= \left\{u\in \mathscr{S}': \sum_{\alpha=0,1} \sum^n_{\ell=1} \|\partial^\alpha_{x_\ell} u\|^{\rm sm \, \cap \, max  \cap \, str} \lesssim \delta \right\}. \label{X1dg}
\end{align}
Moreover,  if $s<1+1/m$ and $F=|\partial_{x_1}u|^m \partial_{x_1}u$, $m\in 2\mathbb{N}$, then \eqref{dnls} is ill-posed in the sense that the solution map $u_0\to u$ is not $C^{m+1}$ from $M^s_{2,1}$ to $C([0,T]; M^s_{2,1})$ for any $T>0$.
\end{thm}

Theorems \ref{th1deri} and \ref{thmderig} also hold in 1D in the cases $\kappa\ge 2$ and $m\ge 4$, respectively for the well posedness with small data in $M^{1/2\kappa}$ and $M^{1+1/m}_{2,1}$  for \eqref{nonl} and \eqref{dnls}; cf. \cite{WaHu07}.
For cubic case in 2D, we need more regularity for the initial data.  First, we need the following semi-norms:
\begin{align}
  \|u\|^{\rm sm2} = & \sum_{k\in \mathbb{Z}^2, \, |k_1|\ge |k_2|\vee 20} \langle k_1\rangle^{3/2} \|\Box_k u\|_{L^\infty_{x_1}L^2_{x_2,t}}  \nonumber \\
  & \quad +  \sum_{k\in \mathbb{Z}^2, \, |k_2|\ge |k_1|\vee 20} \langle k_2\rangle^{3/2} \|\Box_k u\|_{L^\infty_{x_2}L^2_{x_1,t}},  \label{smoothnorm2d}\\
 \|u\|^{\rm max} = &  \sum_{k\in \mathbb{Z}^2}  \left( \|\Box_k u\|_{L^2_{x_1}L^\infty_{x_2,t}}+    \|\Box_k u\|_{L^2_{x_2}L^\infty_{x_1,t}} \right),\\
  \|u\|^{\rm ant} = & \sum_{k\in \mathbb{Z}^2}   \left(\|\Box_k u\|_{L^2_{x_1}L^4_{x_2,t}}+   \|\Box_k u\|_{L^2_{x_2}L^4_{x_1,t}} \right),\\
  \|u\|^{\rm str2} =&  \sum_{k\in \mathbb{Z}^2} \langle k \rangle   \|\Box_k u\|_{L^\infty_{t}L^2_{x} \,  \cap \, L^4_{x,t}},  \ \
  \|u\|^{\rm gstr}  =   \sum_{k\in \mathbb{Z}^2}   \|\Box_k u\|_{ L^3_{x,t}}. \label{gstrnorm2d}
\end{align}
The regularity of $\|\cdot\|^{\rm sm2}$ and $\|u\|^{\rm str2}$ is higher than that of $\|\cdot\|^{\rm sm}$ and  $\|u\|^{\rm str}$.

\begin{thm}\label{thg1} Let $n=2$, $m=2$.
Let $u_0\in M^2_{1,1}$ and there exists a suitably small $\delta>0$ such that $\|u_0\|_{M^2_{1,1}} \le \delta$. Then \eqref{dnls} with nonlinearity \eqref{nonlpoly}  has a unique solution $u \in C(\mathbb{R}, M^2_{2,1}) \cap C_{\rm loc}(\mathbb{R}, M^{3/2}_{1,1}) \cap Y$, where
 \begin{align}
Y= \left\{u\in \mathscr{S}':  \sum_{\alpha=0,1}\sum_{i=1,2} \|\partial^\alpha_{x_i} u\|^{\rm sm2 \, \cap \, max \, \cap \, ant \cap \, str2 \cap \,  gstr} \lesssim \delta \right\}.  \label{X}
\end{align}
In particular, if $u_0\in H^{s,b}$ and $\|u_0\|_{H^{s,b}} \le \delta$ with $s>3$ and $b>1$,  then the result holds.
\end{thm}

\begin{cor}\label{thmmap}
Let $n=2$, $s_0 =(s_1(0), s_2(0), s_3(0))\in \mathbb{S}^2$ with $s_1(0), s_2(0)  \in M^2_{1,1}$ and there exists a suitably small $\delta>0$ such that $\|s_i(0)\|_{M^2_{1,1}} \le \delta$ for $i=1,2$. Then \eqref{Hei} has a unique solution $s=(s_1,s_2,s_3)\in \mathbb{S}^2$ with $s_1,s_2, |s_3|-1 \in C(\mathbb{R}, M^2_{2,1}) \cap C_{\rm loc}(\mathbb{R}, M^{3/2}_{1,1}) \cap X$, where $X$ is as in \eqref{X}.
 \end{cor}

\begin{thm} \label{thm1dm=3}
Let $n=1$, $m=3$. Assume that $u_0 \in M^{11/6}_{1,1}$ and there exists a small $\delta>0$ such that $\|u_0\|_{M^{11/6}_{1,1}} \le \delta$.
Then \eqref{dnls} with nonlinearity \eqref{nonlpoly}  has a unique solution $u \in C(\mathbb{R}, M^{11/6}_{2,1}) \cap C_{\rm loc}(\mathbb{R}, M^{4/3}_{1,1}) \cap Z$, where
 \begin{align}
Z= \left\{u\in \mathscr{S}':  \sum_{\alpha=0,1}  \|\partial^\alpha_{x} u\|^{\rm sm1 \, \cap \, max1 \, \cap \, ant1 \cap \, str1 \cap \,  gstr1} \lesssim \delta \right\}   \label{Z}
\end{align}
and
\begin{align}
 & \|u\|^{\rm sm1} = \sum_{k\in \mathbb{Z}, \, |k|\ge   20} \langle k \rangle^{4/3} \|\Box_k u\|_{L^\infty_{x }L^2_{ t}},  \nonumber \\
 & \|u\|^{\rm max1} = \sum_{k\in \mathbb{Z} }  \|\Box_k u\|_{L^3_{x }L^\infty_{ t}}, \ \
  \|u\|^{\rm ant1} =   \sum_{k\in \mathbb{Z} }    \|\Box_k u\|_{L^3_{x }L^6_{ t}} , \nonumber \\
 &  \|u\|^{\rm str1} =  \sum_{k\in \mathbb{Z} } \langle k \rangle^{5/6}   \|\Box_k u\|_{L^\infty_{t}L^2_{x} \,  \cap \, L^6_{x,t}},  \ \
  \|u\|^{\rm gstr1}  =  \sum_{k\in \mathbb{Z} }   \|\Box_k u\|_{ L^4_{x,t}}.  \nonumber
\end{align}
\end{thm}

\subsection{Strategy of the proof}

We now sketch our ideas in the proof for the main results in the 2D cubic case. We consider the following equivalent integral equation,
 \begin{align}
 u(t)= S(t) u_0 - {\rm i} \int^t_0  S(t-\tau) F(u(\tau))d\tau, \quad S(t):=e^{{\rm i} t\Box}, \ \ \Box=\partial^2_{x_1} - \partial^2_{x_2},
\end{align}
where we assume for simply that $F(u )=\partial_{x_1} (|u|^2 u)$.  By following the smooth effects in 1D as in Kenig-Ponce-Vega \cite{KePoVe93}, the global smooth effects for the solutions of the Schr\"odinger equation in 2D were essentially obtained by Linares and Ponce \cite{LiPo93},
\begin{align}
& \|D^{1/2}_{x_1} S(t) \phi\|_{L^\infty_{x_1}
L^2_{x_2,t} (\mathbb{R}^{1+2})} \lesssim \|\phi\|_{L^2(\mathbb{R}^2)}, \nonumber\\
& \left\|\partial_{x_1} \int^t_0 S(t-\tau) F(\tau) d\tau \right\|_{L^\infty_{x_1} L^2_{x_2,t}  (\mathbb{R}^3)}
\lesssim \|F\|_{L^1_{x_1}  L^2_{x_2,t} (\mathbb{R}^3)}.  \nonumber
\end{align}
The second  estimate from $L^1_{x_1}  L^2_{x_2,t}$ to $ L^\infty_{x_1} L^2_{x_2,t}$ has absorbed one order derivative, which enables us dealing with the derivative in the nonlinearity $F=\partial_{x_1}(|u|^2u)$. By H\"older's inequality,
\begin{align}
  \|  |u|^2 u \|_{L^1_{x_1}L^2_{x_2,t}} \le \|u\|_{L^\infty_{x_1}L^2_{x_2,t}} \|u\|^2_{L^2_{x_1}L^\infty_{x_2,t}}. \label{2dest-1}
\end{align}
So,  one needs at least to estimate $\|u\|_{L^2_{x_1}  L^\infty_{x_2,t}}$. According to the integral equation, we need to show that
$\|S(t)u_0\|_{L^2_{x_1}  L^\infty_{x_2,t}}$ is bounded.  Following \cite{IoKe07}, it was shown in \cite{WaHaHu09} that
\begin{align}
& \|\Box_k S(t) u_0\|_{L_{x_1}^{p}L^{\infty}_{x_2,t} (\mathbb{R}^{1+2} )}   \lesssim \langle k\rangle^{1/p} \|\Box_k u_0\|_{L^2(\mathbb{R}^2)}, \ \ \forall \ p>2. \label{2dest11}
\end{align}
When $p=2$, there is a logarithmic divergence and we can not obtain the time-global estimate. To overcome this difficulty, we will use the Gabor frame. Roughly speaking, any function $u_0$ in $L^2$ and in any modulation space can be expressed in the form $\sum_{k,l\in \mathbb{Z}^2} c_{kl} e^{{\rm i} kx} e^{-|x-l|^2/2}$, it follows that
\begin{align}
S(t) u_0 =  \sum_{k,l\in \mathbb{Z}^2} c_{kl} e^{{\rm i}kx} e^{{\rm i}t (k_1^2- k^2_2)} \prod^2_{j=1} \frac{e^{-\frac{|x_1-l_1 + 2 t  k_1  |^2}{2(1 - 2{\rm i}  t)}   -\frac{|x_2-l_2 - 2 t  k_2  |^2}{2(1 + 2{\rm i}  t)} } }  {(1 + 4 t^2)^{1/2} }.
\end{align}
We can get the following time-global estimates
\begin{align}
& \|\Box_k S(t) u_0\|_{L_{x_1}^{2}L^{\infty}_{x_2,t} (\mathbb{R}^{1+2} )}   \lesssim \langle k\rangle^{1/2} \|\Box_k u_0\|_{L^1(\mathbb{R}^2)}, \label{2dest0}\\
& \|\Box_k \int S(t-\tau)f(\tau) d\tau\|_{L_{x_1}^{2}L^{\infty}_{x_2,t} (\mathbb{R}^{1+2} )}   \lesssim \langle k\rangle^{1/2} \|\Box_k f\|_{L^1_{x,t}(\mathbb{R}^{1+2})}. \label{2dest}
\end{align}
Noticing that  $\|\Box_k u\|_2 \le  \|\Box_k u\|_1$ and comparing \eqref{2dest0} with \eqref{2dest11}, we see that there is  a loss of spatial index in \eqref{2dest0} and \eqref{2dest}. According to \eqref{2dest}, one need to bound (after ignoring the frequency localization and spatial regularity index)
\begin{align}
  \| \partial_{x_1}(|u|^2 u)\|_{L^1_{x,t}(\mathbb{R}^{1+2})}\le \|u_{x_1}\|_{L^\infty_{x_1}L^2_{x_2,t}} \|u\|^2_{L^2_{x_1}L^4_{x_2,t}}. \label{2dest2}
\end{align}
In another way,
\begin{align}
  \| \partial_{x_1}(|u|^2 u)\|_{L^1_{x,t}(\mathbb{R}^{1+2})}\le \|u_{x_1}\|_{L^3_{x ,t}} \|u\|^2_{ L^3_{x ,t}}. \label{2dest3}
\end{align}
\eqref{2dest2} is beneficial to the higher frequency part and \eqref{2dest3} is useful for the lower frequency part. In summary, $L^\infty_{x_1}
L^2_{x_2,t} $ is used for absorbing the derivative in nonlinearity. For the lower frequency part, $L^\infty_{x_1}
L^2_{x_2,t} $ is a bad space and we use the Strichartz space $L^4_{x,t} \cap L^\infty_tL^2_x$ as a substitution.   $L^2_{x_1}
L^\infty_{x_2,t} $ is a maximal function space arising from the nonlinear estimate \eqref{2dest-1}. In order to get a time-global estimate for the Schr\"odinger semi-group in  $L^2_{x_1}
L^\infty_{x_2,t} $,  an intermediate space $L^1_{x,t}$ is
 introduced. Finally,  the anisotropic space $L^2_{x_1}L^4_{x_2,t}$ and the generalized Strichartz space $L^3_{x ,t}$ is employed for the nonlinear estimates in $L^1_{x,t}$.

\section{Linear estimates via Gabor frame }

Gabor frame is of importance in the time-frequency theory, its discrete form enable us to get an exact expression for the solution  of the free Schr\"odinger equation, see below \eqref{discrexp}. Indeed,   Cordero, Nicola  and  Rodino \cite{CoNiRo09} calculated  $e^{{-\rm i} t\Delta} (e^{{\rm i}\alpha kx} e^{- \beta |x-l|^2 })$.  The advantage of the Gabor frame expression  is that it has no singularity at $t=0$ and easier to calculate than the following form
$$
e^{{-\rm i} t\Delta} u_0 = c t^{-n/2} \int e^{ {\rm i} c |x-y|^2/4t} u_0(y) dy.
$$
In this section, we always denote $|\xi|^2_\pm = \sum^n_{j=1} \varepsilon_j \xi^2_j$,  where $\varepsilon_j \in \{1,-1\}$ is arbitrary. For any $x\in \mathbb{R}^n$, we write $\bar{x}=(x_2,...,x_n)$.
\begin{align}
S(t)= \mathscr{F}^{-1}e^{{-\rm i} t|\xi|^2_\pm}
 \mathscr{F},\quad \mathscr{A}  f(t,
x)=\int_0^t S (t-\tau)f(\tau, x)d \tau.
\end{align}

\begin{prop}{\rm (Gabor frame expression)} \label{gfexpress}
Let $s\in \mathbb{R}$, $1\le p,q<\infty$, $u_0 \in M^s_{p,q}$ and
$$
u_0 (x) = \sum_{k,l\in \mathbb{Z}^n} c_{kl} e^{{\rm i}kx} e^{-\frac{|x-l|^2}{2}}.
$$
Then we have
\begin{align} \label{discrexp}
S(t) u_0 =  \sum_{k,l\in \mathbb{Z}^n} c_{kl} e^{{\rm i}kx} e^{{\rm i}t |k|^2_\pm} \prod^n_{j=1} \frac{e^{-\frac{|x_j-l_j +2 t\varepsilon_j k_j  |^2}{2(1-2{\rm i}\varepsilon_j t)}} }{(1-2{\rm i}\varepsilon_j t)^{1/2} }.
\end{align}
\end{prop}
{\bf Proof.} In view of
$$
\widehat{e^{{\rm i}kx} f} = \widehat{f} (\cdot-k), \ \ \widehat{f(\cdot-l)} =e^{-{\rm i}l\xi} \widehat{f}, \ \ \widehat{e^{-|x|^2/2}} = e^{-|\cdot|^2/2},
$$
we see that
$$
\widehat{u}_0 (\xi)= \sum_{k,l\in \mathbb{Z}^n} c_{kl}  e^{{\rm i}lk} e^{-{\rm i}l\xi} e^{-\frac{|\xi-k|^2}{2}}.
$$
It follows that
$$
S(t) u_0 = \sum_{k,l\in \mathbb{Z}^n}  c_{kl}   e^{{\rm i}lk} \mathscr{F}^{-1} \left( e^{{\rm i}t |\xi|^2_\pm} e^{-{\rm i} l \xi} e^{-\frac{|\xi-k|^2}{2}} \right).
$$
In view of
$$
\mathscr{F}^{-1}(f(\cdot-k))   =e^{{\rm i}kx} \mathscr{F}^{-1} f ,   \ \ \mathscr{F}^{-1}(e^{-{\rm i} l \xi} f)  =  (\mathscr{F}^{-1}f) (\cdot-l),
$$
we see that
\begin{align}
S(t) u_0 &  = \sum_{k,l\in \mathbb{Z}^n}  c_{kl}   e^{{\rm i}xk} \left[\mathscr{F}^{-1} \left( e^{{\rm i}t |\xi+k|^2_\pm}   e^{-\frac{|\xi|^2}{2}} \right)\right](x-l) \nonumber\\
&  = \sum_{k,l\in \mathbb{Z}^n}  c_{kl}   e^{{\rm i}xk}  e^{{\rm i}t |k|^2_\pm} \prod^n_{j=1}\left[\mathscr{F}^{-1}_{\xi_j} \left( e^{{\rm i}t \varepsilon_j \xi_j^2 }   e^{-\frac{ \xi_j^2}{2}} \right)\right](x_j-l_j+ 2t \varepsilon_jk_j).
\end{align}
Using the fact that $\mathscr{F}^{-1} e^{-c\xi_1^2/2} = c^{-1/2}e^{-x_1^2/2c}$, we immediately have the result, as desired. $\hfill\Box$

\begin{prop}\label{p1.3} \label{gfexpress0}
Let $n \ge 1$ and $1\le r,p, \bar{p} \leq \infty$. Assume that one of the following alternative conditions holds:
\begin{align}
n\left(\frac{1}{r}- \frac{1}{2} - \frac{1}{\bar{p}}\right) & >\frac{1}{p},  \ \ r\le p, \ \  or \label{rpbarp1}\\
n\left(\frac{1}{r}- \frac{1}{2} - \frac{1}{\bar{p}}\right) & =\frac{1}{p}, \ \   r<p<\infty.  \label{rpbarp2}
\end{align}
Then  we have
\begin{align}
\|S(t) u_0\|_{L_{x_1}^{p}L^{\bar{p}}_{\bar{x},t} (\mathbb{R}^{1+n} )} & \lesssim \|u_0\|_{M^{1/p+ 1- 1/r }_{r,1}}. \label{gabest100}\\
\|\mathscr{A} f\|_{L_{x_1}^{p}L^{\bar{p}}_{\bar{x},t} (\mathbb{R}^{1+n} )} & \lesssim \|f\|_{L^1 (\mathbb{R}, \ M^{1/p+1- 1/r }_{r,1}(\mathbb{R}^n))}. \label{gabest200}
\end{align}
\end{prop}
\noindent{\bf Proof.} Let $u_0 = \sum_{k, l\in \mathbb{Z}^n} c_{kl} e^{{\rm i}kx} e^{-|x-l|^2/2}.$ By Proposition \ref{gfexpress},
\begin{align}
\|S(t)& u_0\|_{L_{x_1}^{p}L^{\bar{p}}_{\bar{x},t}}  \leqslant   \left\| \sum_{k,l\in \mathbb{Z}^n} |c_{kl}|  (1 +|t| )^{-n/2} \prod^n_{j=1}    e^{-\frac{|x_j-l_j +2 t\varepsilon_j k_j  |^2}{2(1+4 t^2)}}\right\|_{L_{x_1}^{p}L^{\bar{p}}_{\bar{x},t}} \nonumber\\
& \leqslant  \sum_{k \in \mathbb{Z}^n}  \left\| \sum_{l_1\in \mathbb{Z}} \langle t \rangle^{-n/2}  e^{-\frac{|x_1-l_1 +2 t\varepsilon_1 k_1  |^2}{2(1+4 t^2)}} \left\| \sum_{\bar{l} \in \mathbb{Z}^{n-1}} |c_{kl}|   \prod^n_{j=2}    e^{-\frac{|x_j-l_j +2 t\varepsilon_j k_j  |^2}{2(1+4 t^2)}}\right\|_{L^{\bar{p}}_{\bar{x}}}\right\|_{L_{x_1}^{p}L^{\bar{p}}_{t}}. \label{barp0}
\end{align}
Applying the fact that $\sup_{x>0} \langle x\rangle^N/e^x <\infty$ for any $N>0$, we have
\begin{align}
  \left\| \sum_{\bar{l} \in \mathbb{Z}^{n-1}} |c_{kl}|   \prod^n_{j=2}    e^{-\frac{|x_j-l_j +2 t\varepsilon_j k_j  |^2}{2(1+4 t^2)}}\right\|_{L^{\bar{p}}_{\bar{x}}} \lesssim   \left\| \sum_{\bar{l} \in \mathbb{Z}^{n-1}} |c_{kl}|  \prod^n_{j=2}  \left\langle  \frac{|x_j -l_j  +2 t \varepsilon_j k_j  |}{ 1+|t|} \right\rangle^{-2} \right\|_{L^{\bar{p}}_{\bar{x}}}.
\end{align}
 In view of  Lemma \ref{convol} we have (see Appendix)
\begin{align} \label{barp1}
  \left\| \sum_{\bar{l} \in \mathbb{Z}^{n-1}} |c_{kl}|   \prod^n_{j=2}    e^{-\frac{|x_j-l_j +2 t\varepsilon_j k_j  |^2}{2(1+4 t^2)}}\right\|_{L^{\bar{p}}_{\bar{x}}} \lesssim \langle t\rangle^{(n-1)/r'+ (n-1)/\bar{p}}  \left(\sum_{\bar{l} \in \mathbb{Z}^{n-1}} |c_{kl}|^r \right)^{1/r}.
\end{align}
It follows from \eqref{barp0} and \eqref{barp1} that
\begin{align} \label{pbarp0}
& \|S(t)  u_0\|_{L_{x_1}^{p}L^{\bar{p}}_{\bar{x},t}} \nonumber\\
& \lesssim  \sum_{k \in \mathbb{Z}^n}  \left\| \sum_{l_1\in \mathbb{Z}} \langle t \rangle^{-\frac{n}{2}+ \frac{n-1}{r'}+ \frac{n-1}{\bar{p}}}  \left(\sum_{\bar{l} \in \mathbb{Z}^{n-1}} |c_{kl}|^r \right)^{1/r}  e^{-\frac{|x_1-l_1 +2 t\varepsilon_1 k_1  |^2}{2(1+4 t^2)}}     \right\|_{L_{x_1}^{p}L^{\bar{p}}_{t}} \nonumber\\
& =  \sum_{k \in \mathbb{Z}^n, \, k_1\neq 0}  \left\| \sum_{l_1\in \mathbb{Z}} \langle t \rangle^{-\frac{n}{2}+ \frac{n-1}{r'}+ \frac{n-1}{\bar{p}}}  \left(\sum_{\bar{l} \in \mathbb{Z}^{n-1}} |c_{kl}|^r \right)^{1/r}  e^{-\frac{|x_1-l_1 +2 t\varepsilon_1 k_1  |^2}{2(1+4 t^2)}}     \right\|_{L_{x_1}^{p}L^{\bar{p}}_{t}} \nonumber\\
& \quad +  \sum_{k \in \mathbb{Z}^n, \, k_1= 0}  \left\| \sum_{l_1\in \mathbb{Z}} \langle t \rangle^{-\frac{n}{2}+ \frac{n-1}{r'}+ \frac{n-1}{\bar{p}}}  \left(\sum_{\bar{l} \in \mathbb{Z}^{n-1}} |c_{kl}|^r \right)^{1/r}  e^{-\frac{|x_1-l_1  |^2}{2(1+4 t^2)}}     \right\|_{L_{x_1}^{p}L^{\bar{p}}_{t}} \nonumber\\
&:= A_{\rm hi} + A_{\rm lo}
\end{align}
We consider the estimate of $A_{\rm hi}$.
\begin{align} \label{pbarp0hi}
 A_{\rm hi} & \lesssim  \sum^3_{s=1}  \sum_{k \in \mathbb{Z}^n, \, k_1\neq 0}  \left\| \sum_{l_1\in \mathbb{Z}} \langle t \rangle^{-\frac{n}{2}+ \frac{n-1}{r'}+ \frac{n-1}{\bar{p}}}  \left(\sum_{\bar{l} \in \mathbb{Z}^{n-1}} |c_{kl}|^r \right)^{1/r}  e^{-\frac{|x_1-l_1 +2 t\varepsilon_1 k_1  |^2}{2(1+4 t^2)}}     \chi_{ \mathbb{D}_s}\right\|_{L_{x_1}^{p}L^{\bar{p}}_{t}} \nonumber\\
& := \Upsilon_1+\Upsilon_2+\Upsilon_3,
\end{align}
where we denote
\begin{align}
\mathbb{D}_1 & = \{(x_1,t): \  |x_1-l_1| > 4|tk_1|\}, \nonumber\\
\mathbb{D}_2 & = \{(x_1,t): \  |x_1-l_1| < |tk_1|\}, \nonumber\\
\mathbb{D}_3 & = \{(x_1,t): \  |tk_1| \leq |x_1-l_1| \leq 4|tk_1|\}. \nonumber
\end{align}
We now estimate $\Upsilon_1$.  Usingt $\sup_{x>0}  \langle x\rangle^\theta/e^x <\infty$ for any $\theta>0$, we have
\begin{align}
\Upsilon_1 &
\lesssim \sum_{k \in \mathbb{Z}^n, \, k_1\neq 0}  \left\| \sum_{l_1\in \mathbb{Z}} \langle t \rangle^{-\frac{n}{2}+ \frac{n-1}{r'}+ \frac{n-1}{\bar{p}}}  \left(\sum_{\bar{l} \in \mathbb{Z}^{n-1}} |c_{kl}|^r \right)^{1/r}  e^{-\frac{|x_1-l_1  |^2}{8(1+4 t^2)}}     \chi_{ \mathbb{D}_1}\right\|_{L_{x_1}^{p}L^{\bar{p}}_{t}} \nonumber\\
& \lesssim  \sum_{k \in \mathbb{Z}^n, \, k_1\neq 0}  \left\| \sum_{l_1\in \mathbb{Z}} \langle x_1-l_1   \rangle^{-\frac{n}{2}+ \frac{n-1}{r'}+ \frac{n-1}{\bar{p}}}  \left(\sum_{\bar{l} \in \mathbb{Z}^{n-1}} |c_{kl}|^r \right)^{1/r} \chi_{ \mathbb{D}_1}\right\|_{L_{x_1}^{p}L^{\bar{p}}_{t}} \nonumber\\
& \lesssim  \sum_{k \in \mathbb{Z}^n, \, k_1\neq 0}  \left\| \sum_{l_1\in \mathbb{Z}} \langle x_1-l_1   \rangle^{-\frac{n}{2}+ \frac{n-1}{r'}+ \frac{n}{\bar{p}}}  \left(\sum_{\bar{l} \in \mathbb{Z}^{n-1}} |c_{kl}|^r \right)^{1/r} \right\|_{L_{x_1}^{p}}. \label{pbarp3}
\end{align}
By Lemma \ref{convol} and \eqref{pbarp3}, we have
$$
\Upsilon_1 \lesssim  \sum_{k \in \mathbb{Z}^n}  \left(\sum_{ l  \in \mathbb{Z}^{n }} |c_{kl}|^r \right)^{1/r}.
$$
Noticing that $|t| \sim |x_1-l_1|/|k_1|$ in $\mathbb{D}_3$, we have
\begin{align}
\Upsilon_3
 & \lesssim  \sum_{k \in \mathbb{Z}^n, \, k_1\neq 0}  \left\| \sum_{l_1\in \mathbb{Z}} \left\langle \frac{ x_1-l_1 }{|k_1|}    \right\rangle^{-\frac{n}{2}+ \frac{n-1}{r'}+ \frac{n-1}{\bar{p}}}  \left(\sum_{\bar{l} \in \mathbb{Z}^{n-1}} |c_{kl}|^r \right)^{1/r} \chi_{ \mathbb{D}_3}\right\|_{L_{x_1}^{p}L^{\bar{p}}_{t}} \nonumber\\
& \lesssim  \sum_{k \in \mathbb{Z}^n, \, k_1\neq 0}  \left\| \sum_{l_1\in \mathbb{Z}}  \left\langle \frac{ x_1-l_1 }{|k_1|}    \right\rangle^{-\frac{n}{2}+ \frac{n-1}{r'}+ \frac{n}{\bar{p}}}   \left(\sum_{\bar{l} \in \mathbb{Z}^{n-1}} |c_{kl}|^r \right)^{1/r} \right\|_{L_{x_1}^{p}}. \label{pbarp4}
\end{align}
By Lemma \ref{convol}, we have
\begin{align}
\Upsilon_3
 & \lesssim    \sum_{k \in \mathbb{Z}^n}  \langle  k_1 \rangle^{  \frac{1}{r'}+ \frac{1}{ p }}   \left(\sum_{ l  \in \mathbb{Z}^{n }} |c_{kl}|^r \right)^{1/r}  . \label{pbarp5}
\end{align}
For the estimate of $\Upsilon_2$, we have
\begin{align}
\Upsilon_2 &
\lesssim \sum_{k \in \mathbb{Z}^n, \, k_1\neq 0}  \left\| \sum_{l_1\in \mathbb{Z}} \langle t \rangle^{-\frac{n}{2}+ \frac{n-1}{r'}+ \frac{n-1}{\bar{p}}}  \left(\sum_{\bar{l} \in \mathbb{Z}^{n-1}} |c_{kl}|^r \right)^{1/r}  \chi_{ \mathbb{D}_2}\right\|_{L_{x_1}^{p}L^{\bar{p}}_{t}} \nonumber\\
&
\lesssim \sum_{k \in \mathbb{Z}^n, \, k_1\neq 0}  \left\| \sum_{l_1\in \mathbb{Z}} \left\langle \frac{ x_1-l_1 }{|k_1|}  \right\rangle^{-\frac{n}{2}+ \frac{n-1}{r'}+ \frac{n}{\bar{p}}}  \left(\sum_{\bar{l} \in \mathbb{Z}^{n-1}} |c_{kl}|^r \right)^{1/r}  \right\|_{L_{x_1}^{p}}. \label{pbarp6}
\end{align}
By Lemma \ref{convol}, we see that $\Upsilon_2$ has the same upper bound as $\Upsilon_3$ in \eqref{pbarp5}.  Collecting the estimates of $\Upsilon_1$, $\Upsilon_2$ and $\Upsilon_3$, we have the desired estimate.

We consider the estimate of $A_{\rm lo}$.  We have
\begin{align} \label{pbarp0}
 A_{\rm lo}
& \lesssim  \sum^2_{s=1}  \sum_{k \in \mathbb{Z}^n, \, k_1= 0}  \left\| \sum_{l_1\in \mathbb{Z}} \langle t \rangle^{-\frac{n}{2}+ \frac{n-1}{r'}+ \frac{n-1}{\bar{p}}}  \left(\sum_{\bar{l} \in \mathbb{Z}^{n-1}} |c_{kl}|^r \right)^{1/r}  e^{-\frac{|x_1-l_1 |^2}{2(1+4 t^2)}}     \chi_{ \mathbb{E}_s}\right\|_{L_{x_1}^{p}L^{\bar{p}}_{t}} \nonumber\\
& := \Xi_1+\Xi_2,
\end{align}
where we denote
\begin{align}
\mathbb{E}_1 & = \{(x_1,t): \  |x_1-l_1| > |t|\}, \nonumber\\
\mathbb{E}_2 & = \{(x_1,t): \  |x_1-l_1| \leq |t |\}. \nonumber
\end{align}
Applying the fact $\sup_{x>0} \langle x\rangle^\theta/e^x <\infty$, we have
\begin{align} \label{pbarp00}
 \Xi_1  & \lesssim    \sum_{k \in \mathbb{Z}^n, \, k_1= 0}  \left\| \sum_{l_1\in \mathbb{Z}} \langle x_1-l_1 \rangle^{-\frac{n}{2}+ \frac{n-1}{r'}+ \frac{n-1}{\bar{p}}}  \left(\sum_{\bar{l} \in \mathbb{Z}^{n-1}} |c_{kl}|^r \right)^{1/r}       \chi_{ \mathbb{E}_1}\right\|_{L_{x_1}^{p}L^{\bar{p}}_{t}} \nonumber\\
&  \lesssim    \sum_{k \in \mathbb{Z}^n, \, k_1= 0}  \left\| \sum_{l_1\in \mathbb{Z}} \langle x_1-l_1 \rangle^{-\frac{n}{2}+ \frac{n-1}{r'}+ \frac{n}{\bar{p}}}  \left(\sum_{\bar{l} \in \mathbb{Z}^{n-1}} |c_{kl}|^r \right)^{1/r}       \right\|_{L_{x_1}^{p} }.
\end{align}
By Lemma \ref{convol}, we have
$$
\Xi_1 \lesssim  \sum_{\bar{k} \in \mathbb{Z}^{n-1}}  \left(\sum_{ l  \in \mathbb{Z}^{n }} |c_{(0,\bar{k})l}|^r \right)^{1/r}.
$$
Analogous to the estimate of $\Upsilon_2$, we can show that $\Xi_2$ has the same upper bound as $\Xi_1$. $\hfill\Box$

\medskip

Let us observe an endpoint case $r=1$. We have

\begin{cor}  \label{gfexpress1}
Let $n \ge 1$, $1\le p, \  \bar{p} \le \infty$. Assume one of the following alternative conditions holds:
$$
  n\left(\frac{1}{2} - \frac{1}{\bar{p}}\right) > \frac{1}{p}; \ \ or
$$
$$
  n\left(\frac{1}{2} - \frac{1}{\bar{p}}\right) = \frac{1}{p},  \ \  1<p<\infty.
$$
Then  we have
\begin{align}
\|S(t) u_0\|_{L_{x_1}^{p}L^{\bar{p}}_{\bar{x},t} (\mathbb{R}^{1+n} )} & \lesssim \|u_0\|_{M^{1/p}_{1,1}}. \label{gabest1}\\
\|\mathscr{A} f\|_{L_{x_1}^{p}L^{\bar{p}}_{\bar{x},t} (\mathbb{R}^{1+n} )} & \lesssim \|f\|_{L^1 (\mathbb{R}, \ M^{1/p}_{1,1}(\mathbb{R}^n))}. \label{gabest2}
\end{align}
\end{cor}

\section{Linear estimates with $\Box_k$-decomposition}

\begin{cor}\label{loc-gfexpress}
{\rm ($L^1$-anisotropic estimates)} Let $n \ge 1$, $1\le p, \bar{p}\le \infty$.
Assume one of the following alternative conditions holds:
$$
  n\left(\frac{1}{2} - \frac{1}{\bar{p}}\right) > \frac{1}{p}; \ \ or
$$
$$
  n\left(\frac{1}{2} - \frac{1}{\bar{p}}\right) = \frac{1}{p},  \ \  1<p<\infty.
$$
Then  we have
\begin{align}
\|\Box_k S(t) u_0\|_{L_{x_1}^{p}L^{\bar{p}}_{\bar{x},t} (\mathbb{R}^{1+n} )} & \lesssim \langle k\rangle^{1/p} \|\Box_k u_0\|_{L^1(\mathbb{R}^n)},  \label{fre-local-1}\\
\|\Box_k \mathscr{A} f\|_{L_{x_1}^{p}L^{\bar{p}}_{\bar{x},t} (\mathbb{R}^{1+n} )} & \lesssim \langle k\rangle^{1/p} \|\Box_k f\|_{L^1_{x,t} (\mathbb{R}^{n+1})}. \label{fre-local-2}
\end{align}
In particular, we have for $n=2$,
\begin{align}
\max_{q=4,\infty}\|\Box_k S(t) u_0\|_{L_{x_1}^{2}L^{q}_{x_2,t} (\mathbb{R}^{1+2} )} & \lesssim \langle k\rangle^{1/2} \|\Box_k u_0\|_{L^1(\mathbb{R}^2)},  \label{fre-local-12d}\\
\max_{q=4,\infty}\|\Box_k \mathscr{A} f\|_{L_{x_1}^{2}L^q_{x_2,t} (\mathbb{R}^{1+2} )} & \lesssim \langle k\rangle^{1/2} \|\Box_k f\|_{L^1_{x,t} (\mathbb{R}^{n+1})}. \label{fre-local-22d}
\end{align}
\begin{align}
\|\Box_k S(t) u_0\|_{L^{3}_{x ,t} (\mathbb{R}^{1+2} )} & \lesssim \langle k\rangle^{1/3} \|\Box_k u_0\|_{L^1(\mathbb{R}^2)},  \label{fre-local-13}\\
\|\Box_k \mathscr{A} f\|_{ L^3_{x,t} (\mathbb{R}^{1+2} )} & \lesssim \langle k\rangle^{1/3} \|\Box_k f\|_{L^1_{x,t} (\mathbb{R}^{n+1})}. \label{fre-local-23}
\end{align}

\end{cor}
\noindent{\bf{Proof.}}
 By Proposition \ref{gfexpress1},  we have
\begin{align}
\|\Box_k S(t) u_0\|_{L_{x_1}^{p}L^{\bar{p}}_{\bar{x},t} (\mathbb{R}^{1+n} )} \lesssim \|\Box_k u_0\|_{M^{1/p}_{1,1}}.
\end{align}
By definition and $\Box_k: L^r \to L^r$\footnote{$|l|_\infty= \max (|l_1|,..., |l_n|)$ for $l=(l_1,..., l_n) \in \mathbb{Z}^n$.},
$$
\|\Box_k u_0\|_{M^{1/p}_{1,1}} \le \sum_{|l|_\infty \leq 1} \langle k+l \rangle^{1/p}\|\Box_{k+l} \Box_k u_0\|_1 \lesssim  \langle k\rangle^{1/p} \|\Box_k u_0\|_1,
$$
which implies the result, as desired.  $\hfill\Box$

\begin{prop} \label{smeff} {\rm (Smooth effects, \cite{LiPo93, WaHaHu09})}
For any  $k=(k_1, \ldots, k_n)\in\mathbb{Z}^n$,  we have
\begin{align}
\|D^{1/2}_{x_1} \Box_k S(t) u_0
 \|_{L_{x_1}^{\infty} L^2_{\bar{x},t} } & \lesssim
 \|\Box_k u_0\|_2,  \label{p1.4}\\
\|\partial_{x_1} \Box_k \mathscr{A}
f\|_{L_{x_1}^{\infty} L^2_{\bar{x},t} } & \lesssim
 \|\Box_k f\|_{L_{x_1}^{1}L_{\bar{x}, t}^2 }.\label{p1.5}
\end{align}
\end{prop}

\begin{prop}\label{strichartz}
{\rm (Strichartz estimates, \cite{WaHe07, WaHaHu09})} Let $4/n \le p < \infty$.  We have
\begin{align}
 & \left \|\Box_k S(t) u_0 \right\|_{
L^{2+p}_{t,x}\, \cap \, L^\infty_t L^2_x (\mathbb{R}^{1+n}) }
\lesssim \|\Box_k u_0\|_{L^2 (\mathbb{R}^{n}) }, \label{st-sm-m-c1}\\
& \left \|\Box_k \mathscr{A}  f \right\|_{L^\infty_t L^2_x \, \cap
\, L^{2+p}_{t,x}  (\mathbb{R}^{1+n}) } \lesssim    \|\Box_k f\|_{
L^{(2+p)/(1+p)}_{t,x} (\mathbb{R}^{1+n}) }. \label{st-sm-m-c9}
\end{align}
\end{prop}

\begin{prop}\label{str-sm-relation}
{\rm (Interaction estimates, \cite{WaHaHu09})} Let $4/n \le p < \infty$.  We have
\begin{itemize}
     \item[\rm (1)] Smooth-Strichartz estimates
\begin{align}
& \left \|\Box_k \partial_{x_1} \mathscr{A}  f \right\|_{L^\infty_t L^2_x \, \cap
\, L^{2+p}_{t,x}  (\mathbb{R}^{1+n}) } \lesssim  \langle
k_1\rangle^{1/2}  \|\Box_k f\|_{ L^1_{x_1} L^2_{\bar{x},t }
(\mathbb{R}^{1+n}) },  \nonumber
\end{align}
\item[\rm (2)] Strichartz-smooth estimates
\begin{align}
&\left \|\Box_k \mathscr{A} \partial_{x_1} f
\right\|_{L^\infty_{x_1} L^2_{\bar{x},t}  (\mathbb{R}^{1+n})}
\lesssim \langle k_1\rangle^{1/2}\| \Box_k f\|_{
L^{(2+p)/(1+p)}_{t,x} (\mathbb{R}^{1+n})}, \nonumber
\end{align}
 \item[\rm (3)] Strichartz-maximal estimates
 \begin{align}
&  \|\Box_k \partial_{x_i} \mathscr{A} f\|_{L^q_{x_1} L^\infty_{\bar{x}, t}} \lesssim \langle k_i\rangle \langle k_1\rangle^{1/q} \|\Box_k f\|_{L^{ (2+p)/(1+p)}_{x,t}}, \ \ 2\le q\le \infty, \nonumber
\end{align}
\item[\rm (4)] Smooth-maximal estimates
\begin{align}
 & \|\Box_k \partial_{x_i} \mathscr{A} f\|_{L^q_{x_1} L^\infty_{\bar{x}, t}} \lesssim \langle k_i\rangle^{1/2} \langle k_1\rangle^{1/q} \|\Box_k f\|_{L^1_{x_i} L^2_{(x_j)_{j\neq i},t}},  \ \ 2<q\le \infty. \nonumber
 \end{align}
\end{itemize}
\end{prop}

\begin{prop}  \label{gfexpress2}
{\rm ($L^2$-anisotropic estimates)} Let $n \ge 1$, $2 \le q, \  \bar{q} \le \infty$. Assume one of the following alternative conditions holds:
$$
  n\left(\frac{1}{2} - \frac{2}{\bar{q}}\right) > \frac{2}{q}; \ \ or
$$
$$
  n\left(\frac{1}{2} - \frac{2}{\bar{q}}\right) = \frac{2}{q},  \ \  2<q<\infty.
$$
Then  we have
\begin{align}
\|\Box_k S(t) u_0\|_{L_{x_1}^{q}L^{\bar{q}}_{\bar{x},t} (\mathbb{R}^{1+n} )} & \lesssim \langle k_1\rangle^{1/q} \|\Box_k u_0\|_{2}, \label{2gabest1}\\
\|\Box_k \mathscr{A} f\|_{L_{x_1}^{q}L^{\bar{q}}_{\bar{x},t} (\mathbb{R}^{1+n} )} & \lesssim \langle k\rangle^{1/q} \|\Box_k f\|_{L^1_tL^{2}_x (  (\mathbb{R}^{n+1})}. \label{2gabest2}
\end{align}
\end{prop}
{\bf Proof.} By duality, it suffices to show that
 \begin{align}
\int_{\mathbb{R}} (\Box_k S(t) u_0, \phi(t)) dt  \lesssim \langle k_1\rangle^{1/q}\|\Box_k u_0\|_{2} \|\phi\|_{L^{q'}_{x_1} L^{\bar{q}'}_{\bar{x},t}}.  \end{align}
It is easy to see that
 \begin{align}
\int_{\mathbb{R}} (\Box_k S(t) u_0, \phi(t)) dt  \lesssim \langle k_1\rangle^{1/q}\|\Box_k u_0\|_{2} \sum_{|l|_\infty \le 1} \left\|\Box_{k+l}\int S(-t) \phi(t) dt \right\|_{2}.
\end{align}
We have
 \begin{align}
  \left\|\Box_k \int S(-t) \phi(t) dt \right\|^2_{2}  \leq  \|\phi\|_{L^{q'}_{x_1} L^{\bar{q}'}_{\bar{x},t}} \left\|\Box_k \int S(t-s) \phi(s) ds \right\|_{L^{q}_{x_1} L^{\bar{q}}_{\bar{x},t}}.
\end{align}
Let us observe that
$$
\Box_k \int S(t-s) \phi(s) ds  = (\mathscr{F}^{-1} e^{{\rm i}t |\xi|^2_\pm} \eta_k (\xi)) * \phi,
$$
where $*$ denotes the convolution on $x$ and $t$.   Applying Young's inequality, we obtain that
 \begin{align}
\left\|\Box_k \int S(t-s) \phi(s) ds \right\|_{L^{q}_{x_1} L^{\bar{q}}_{\bar{x},t}} \lesssim  \left \|\mathscr{F}^{-1} e^{{\rm i}t |\xi|^2_\pm} \eta_k (\xi) \right\|_{L^{q/2}_{x_1} L^{\bar{q}/2}_{\bar{x},t}}  \|\phi\|_{L^{q'}_{x_1} L^{\bar{q}'}_{\bar{x},t}}.
\end{align}
Hence, if we can show that
 \begin{align}
   \|\mathscr{F}^{-1} e^{{\rm i}t |\xi|^2_\pm} \sigma_k (\xi)\|_{L^{q/2}_{x_1} L^{\bar{q}/2}_{\bar{x},t}}  \lesssim \langle k_1\rangle^{2/q},
\end{align}
then the result follows. Indeed, in view of Corollary \ref{gfexpress1},
 \begin{align} \nonumber
   \|\mathscr{F}^{-1} e^{{\rm i}t |\xi|^2_\pm} \sigma_k (\xi)|_{L^{q/2}_{x_1} L^{\bar{q}/2}_{\bar{x},t}}  \lesssim \langle k_1\rangle^{2/q} \|\mathscr{F}^{-1}\sigma_k \|_{1} \lesssim \langle k_1\rangle^{2/q}.   \ \ \ \ \ \ \  \Box
\end{align}
Following some idea  as in Bejenaru,  Ionescu,  Kenig and Tataru's  \cite{BeIoKeTa08}, we have
\begin{prop}  \label{maxsmest}
Suppose that the conditions of Proposition \ref{gfexpress2} are satisfied.
Then for $|k_1|\ge 20$,  we have
\begin{align}
\|\Box_k  \mathscr{A} f\|_{L_{x_1}^{q}L^{\bar{q}}_{\bar{x},t} (\mathbb{R}^{1+n})} & \lesssim \langle k_1\rangle^{1/q-1/2} \|\Box_k f\|_{L^1_{x_1}L^{2}_{\bar{x},t} (  (\mathbb{R}^{n+1})}. \label{maxsm}
\end{align}
\end{prop}
{\bf Proof.}
Denote
\begin{align}
u  = c \mathscr{F}^{-1}_{\tau, \xi} \frac{1}{|\xi|^2_\pm -\tau}
\mathscr{F}_{t,x} f.  \label{maximal1}
\end{align}
 We can assume that $|\xi|^2_\pm =
\xi^2_1+ \varepsilon_2 \xi^2_2+...+ \varepsilon_n \xi^n_n: =\xi^2_1+
|\bar{\xi}|^2_\pm$. Denote
\begin{align}
\mathbb{E} = \{(\tau, \bar{\xi})\in \mathbb{R}^n: \ |\bar{\xi}|^2_\pm - \tau <0\}. \label{E}
\end{align}
Let $\sigma_k(\xi)=\eta_{k_1}(\xi_1)...\eta_{k_n}(\xi_n)$ be as in \eqref{denote} and \eqref{denote2}. We can assume that $k_1>0$. We have
\begin{align}
\Box_k u  &  = c \int_{\mathbb{R}^{1+n}} \frac{e^{{\rm i}t\tau} e^{{\rm i}x \xi} }{|\xi|^2_\pm -\tau}
  \widehat{\Box_k f} (\tau, \xi)  d\xi d\tau \nonumber\\
 &  = c \int_{\mathbb{R}\times \mathbb{E}} \sum_{|\ell_1|\le 10} \eta_{ k_1 +\ell_1} \left(\sqrt{\tau-|\bar{\xi}|^2_\pm} \right)  \frac{e^{{\rm i}t\tau} e^{{\rm i}x \xi} }{|\xi|^2_\pm -\tau}
  \widehat{\Box_k f} (\tau, \xi) d\xi d\tau  \nonumber\\
  & \quad \quad + c \int_{\mathbb{R}\times \mathbb{E}}  \sum_{|\ell_1|> 10} \eta_{ k_1 +\ell_1} \left(\sqrt{\tau-|\bar{\xi}|^2_\pm} \right)    \frac{e^{{\rm i}t\tau} e^{{\rm i}x \xi} }{|\xi|^2_\pm -\tau}\widehat{\Box_k f} (\tau, \xi) d\xi d\tau  \nonumber\\
  & \quad \quad + c \int_{\mathbb{R}\times (\mathbb{R}^n\setminus \mathbb{E})}   \frac{e^{{\rm i}t\tau} e^{{\rm i}x \xi} }{|\xi|^2_\pm -\tau}\widehat{\Box_k f} (\tau, \xi) d\xi d\tau  \nonumber\\
  &:= I+II+III.
\label{maximal2}
\end{align}
Let  $\widehat{f(y_1)} (\tau, \bar{\xi})$ be the Fourier transform of $f(t,y_1, \bar{y})$ with respect to $t$ and $ \bar{y}$.  Using the almost orthogonal property of $\Box_k$, we have
\begin{align}
I  = c \sum_{|\bar{\ell}|_\infty \le 1, \ |\ell_1|\le 10} \int_{\mathbb{R}}  \int_{\mathbb{E}} e^{{\rm i}(t\tau +\bar{x} \bar{\xi})}   \sigma_{\bar{k}+\bar{\ell}} (\bar{\xi}) & \eta_{k_1+\ell_1} \left(\sqrt{\tau-|\bar{\xi}|^2_\pm} \right) \widehat{\Box_k f( y_1)} (\tau, \bar{\xi}) \nonumber\\
 & \times \left(\int_{\mathbb{R}} \frac{e^{{\rm i} (x_1-y_1)\xi_1} }{\xi^2_1+ |\bar{\xi}|^2_\pm -\tau} d\xi_1 \right)
 d\bar{\xi} d\tau dy_1.
  \label{maximal3}
\end{align}
For any $a>0$, one has that
\begin{align}
 \int_{\mathbb{R}} \frac{e^{{\rm i} x_1 \xi_1} }{\xi^2_1-a^2} d\xi_1
 = \frac{{\rm sgn} (x_1)}{2a} (e^{{\rm i} x_1 a}- e^{-{\rm i} x_1 a}).
  \label{maximal4}
\end{align}
 Applying \eqref{maximal4} and changing the variable, the norm in $L^q_{x_1} L^{\bar{q}}_{\bar{x},t}$ of the right hand side of \eqref{maximal3} can be reduced to the following estimate
\begin{align}
\Gamma & := \left\| \int_{\mathbb{R}} {\rm sgn} (x_1-y_1)   \left( \mathscr{F}^{-1}  e^{{\rm i} t |\xi|^2_\pm }    \sigma_{\bar{k} } (\bar{\xi})   \eta_{k_1} \left( \xi_1 \right) \widehat{\Box_k f( y_1)} ( |\xi|^2_\pm, \bar{\xi}) \right)(x_1-y_1, \bar{x})  dy_1 \right\|_{L^q_{x_1} L^{\bar{q}}_{\bar{x},t}} \nonumber\\
  & \lesssim \int_{\mathbb{R}}  \left\|\left( \mathscr{F}^{-1}  e^{{\rm i} t |\xi|^2_\pm }    \sigma_{\bar{k} } (\bar{\xi})   \eta_{k_1} \left( \xi_1 \right) \widehat{\Box_k f( y_1)} ( |\xi|^2_\pm, \bar{\xi}) \right)(x_1-y_1, \bar{x})  \right\|_{L^q_{x_1} L^{\bar{q}}_{\bar{x},t}} dy_1.
  \label{maximal6}
\end{align}
Applying  Proposition \ref{gfexpress2}, we have
 \begin{align}
 \Gamma & \lesssim \int_{\mathbb{R}} \langle k_1\rangle^{1/q}
 \left\| \chi_{[k_1-1, k_1+1]}  ( \xi_1) \widehat{\Box_k f( y_1)} ( |\xi|^2_\pm, \bar{\xi})    \right\|_{L^2} dy_1 \nonumber\\
  & \lesssim  \langle k_1\rangle^{1/q-1/2} \int_{\mathbb{R}}
 \left\| (\Box_k f )(t, y_1, \bar{y})  \right\|_{L^2_{\bar{y},t}} dy_1.
\end{align}
Next, we consider the estimate of $II$. We have
\begin{align}
II  = c \sum_{ |\ell|_{\infty} \le 1, \ |j_1|> 10} \int_{\mathbb{R}}  \int_{\mathbb{E}} e^{{\rm i}(t\tau +\bar{x} \bar{\xi})}  &  \sigma_{\bar{k}+\bar{\ell}} (\bar{\xi}) \eta_{k_1+j_1} \left(\sqrt{\tau-|\bar{\xi}|^2_\pm} \right) \widehat{\Box_k f( y_1)} (\tau, \bar{\xi}) \nonumber\\
 & \times \left(\int_{\mathbb{R}} \frac{e^{{\rm i} (x_1-y_1)\xi_1} \eta_{k_1+\ell_1}(\xi_1) }{\xi^2_1+ |\bar{\xi}|^2_\pm -\tau} d\xi_1 \right)
 d\bar{\xi} d\tau dy_1.
  \label{maximal8}
\end{align}
Integrating by part, we see that
\begin{align}
\left|\int_{\mathbb{R}} \frac{e^{{\rm i} x_1 \xi_1} \eta_{k_1+\ell_1}(\xi_1) }{\xi^2_1-s} d\xi_1 \right| \lesssim \frac{1}{1+|x_1|} \max_{|\xi_1-k_1| \le 1}
 \frac{1}{|\xi^2_1-s|}: =K(x_1,s).  \label{maximal9}
\end{align}
Using $\|\mathscr{F}f\|_{\bar{q}} \lesssim \|f\|_{\bar{q}'}$,  one has that
\begin{align}
& \|II\|_{L^{q}_{x_1} L^{\bar{q}}_{\bar{x},t}} \nonumber\\
& \lesssim   \sum_{ |\ell|_{\infty} \le 1} \int_{\mathbb{R}}   \left \|  \sigma_{\bar{k}+\bar{\ell}} \sum_{|j_1|> 10} \eta_{k_1+j_1}  (\sqrt{s}) \widehat{\Box_k f( y_1)} (s+|\bar{\xi}|^2_\pm, \bar{\xi}) K(x_1-y_1,s)
   \right\|_{L^{q}_{x_1} L^{\bar{q}'}_{\bar{\xi},s>0}}
\!\!\!  dy_1 \nonumber\\
 & \lesssim   \sum_{ |\ell|_{\infty} \le 1} \int_{\mathbb{R}}   \left \| \sigma_{\bar{k}+\bar{\ell}}  \sum_{|j_1|> 10} \eta_{k_1+j_1}  (\sqrt{s}) \widehat{\Box_k f( y_1)} (s+|\bar{\xi}|^2_\pm, \bar{\xi}) [(|k_1|\pm 1)^2-s]^{-1}
   \right\|_{ L^{\bar{q}'}_{\bar{\xi}, s>0}}
\!\!\!  dy_1 \nonumber\\
& \lesssim   \langle k_1\rangle^{-1/2-1/q}  \int_{\mathbb{R}}   \left \| \Box_k f
   \right\|_{ L^{2}_{\bar{x}, t}}
   dy_1.
  \label{maximal10}
\end{align}
Since
$$
\mathscr{A} f= u -c \int_{\mathbb{R}} S(t-s) {\rm sgn}(s) f(s) ds,
$$
and the dual estimate argument implies that (see \cite{WaHaHu09})
$$
\left\|\Box_k \int_{\mathbb{R}} S(t-s) {\rm sgn}(s) f(s) ds \right\|_{L^{q}_{x_1} L^{\bar{q}}_{\bar{x},t}} \lesssim  \langle k_1\rangle^{-1/2 + 1/q}  \|\Box_k f\|_{L^{1}_{x_1} L^{2}_{\bar{x},t} }.
$$
The result follows. $\hfill \Box$

\medskip

Let us recall that for $q>2$, Proposition \ref{maxsmest} was obtained in  \cite{WaHaHu09} by using the standard dual estimate argument. However, in order to get the sharp global well posedness result of \eqref{dnls} with the nonlinearity \eqref{nonlpoly}, the result in the case $q=2$ is of importance.

\begin{prop}\label{ubinL1}
Let $1\le p \leq \infty$. Then $\Box_k S (t) : L^p \to L^p$ is uniformly bounded. More precisely,
\begin{align}
\|\Box_kS (t) u_0\|_{L^p }
 \lesssim (1+|t|^{n/2}) \|\Box_k  u_0\|_{L^p }
\label{modppa11}
\end{align}
uniformly holds for all $k\in \mathbb{Z}^n$.
\end{prop}
{\bf Proof.} See \cite{BeGrOkRo06, WaZhGu06}.

\section{Proof of Theorems \ref{th1deri} and \ref{thmderig}}

 Proposition \ref{maxsmest} is crucial for us to  reach the critical space $M^{1/2\kappa}_{2,1}$ and  $M^{1+1/m}_{2,1}$ for $m=2\kappa=2$ in Theorems \ref{th1deri} and \ref{thmderig}, respectively.   For convenience, we write
\begin{align}
  \|u\|^{\rm sm}_i = &   \sum_{k\in \mathbb{Z}^n, \, |k_i|\ge 20 \vee \max_{j\neq i}|k_j|} \langle k_i\rangle^{1/2+1/m} \|\Box_k u\|_{L^\infty_{x_i}L^2_{(x_j)_{j\neq i},t}} , \\
 \|u\|^{\rm max}_i  = &    \sum_{k\in \mathbb{Z}^n}  \|\Box_k u\|_{L^m_{x_i}L^\infty_{(x_j)_{j\neq i},t}} , \\
  \|u\|^{\rm str} =&  \sum_{k\in \mathbb{Z}^2} \langle k \rangle^{1/m}   \|\Box_k u\|_{L^\infty_{t}L^2_{x} \,  \cap \, L^{m+2}_{x,t}}.
\end{align}
For simplicity, we assume that $\mu=0$.  When $\mu \neq 0$, $|u|^{2\nu} u$ can be handled by only using the Strichartz space.  We consider the mapping
 \begin{align}
\mathcal{T}:  u(t) \to  S(t) u_0 - {\rm i}  \mathscr{A} \vec{\lambda} \cdot \nabla (|u|^{2\kappa} u).
\end{align}

\begin{lem}\label{connorm}
We have for any $k=(k_1,...,k_n)$ with $|k_1| \ge |k_2| \vee 20$ and $1\le p,q \le \infty$,
\begin{align} \label{21con}
\|\partial^\alpha_{x_2}\Box_k f\|_{L^{q}_{x_1}L^p_{\bar{x} ,t}}  \lesssim  \|\partial_{x_1}\Box_k f\|_{L^q_{x_1}L^p_{\bar{x} ,t}}, \ \ \alpha=0,1.
\end{align}
In particular, one has that $\| \partial^\alpha_{x_2} u \|^{\rm sm}_1 \lesssim  \| \partial_{x_1} u \|^{\rm sm}_1$, $\alpha=0,1$.
\end{lem}
{\bf Proof.}  By the almost orthogonality of $\Box_k$ and noticing that $|k_1| \ge |k_2| \vee 20$, we have
\begin{align}
      \|\partial_{x_2}\Box_k f\|_{L^q_{x_1}L^p_{\bar{x} ,t}}
& \lesssim   \sum_{|l|_\infty \le 1} \left\| \mathscr{F}^{-1} \left(\frac{\xi_2}{\xi_1} \sigma_{k+l} \right)\right\|_{1}\|\partial_{x_1}\Box_k u\|_{L^q_{x_1}L^p_{\bar{x} ,t}} \nonumber\\
& \lesssim \|\partial_{x_1}\Box_k u\|_{L^q_{x_1}L^p_{x_2,t}},
\end{align}
which implies the result, as desired.  $\hfill\Box$

\begin{lem}\label{boxkest}
For any $k, k^{(s)} \in \mathbb{Z}^n$, $k^{(s)}=(k^{(s)}_1,..., k^{(s)}_n)$, we have
$  \Box_k \bar{u} = (-1)^n \overline{\Box_{-k} u},$ and
$$
 \Box_k(\Box_{k^{(1)}}u_1...\Box_{k^{(r)}}u_r) =0
$$
if $|k_i-k^{(1)}_i-...-k^{(r)}_i | > r+1$, $i=1, ..., n$.
\end{lem}
{\bf Proof.} See \cite{WaHe07}.

\begin{lem}\label{smthnorm}
Let $m=2\kappa$ with $\kappa \in \mathbb{N}$.   We have
\begin{align} \label{sm11}
\| \mathcal{T} u\|^{\rm sm}_{i}  \lesssim  \|u_0\|_{M^{1/2\kappa}_{2,1}} + \| u\|^{\rm sm} (\| u\|^{\rm max})^{2\kappa} + (\| u\|^{\rm str})^{2\kappa+1}
\end{align}
\end{lem}
{\bf Proof.}  It suffices to bound $\| \mathcal{T} u\|^{\rm sm}_{1} $.  Applying the 1/2-order smoothness of $S(t)$ and Lemma \ref{connorm},
 \begin{align} \label{sm11}
\| \mathcal{T} u\|^{\rm sm}_{1}  \lesssim  \|u_0\|_{M^{1/2\kappa}_{2,1}} + \|\partial_{x_1}\mathscr{A} (|u|^{2\kappa} u)\|^{\rm sm}_1.
\end{align}
For convenience, we write
\begin{align}
& \mathbb{A}^{\lambda, i}_{\rm lo} = \left\{(k^{(1)},..., k^{(\lambda)}) \in (\mathbb{Z}^n)^\lambda: \, \max_{1\le s\le \lambda} |k^{(s)}_i| < 20 \right\}, \nonumber\\
& \mathbb{A}^{\lambda, i}_{\rm hi} = \left\{(k^{(1)},..., k^{(\lambda)}) \in (\mathbb{Z}^n)^\lambda: \, \max_{1\le s\le \lambda} |k^{(s)}_i| \ge 20 \right\}
\end{align}
for $k^{(s)}= ( k^{(s)}_1,...,  k^{(s)}_n)$.  Since  $  \Box_k \bar{u} = \ \overline{\Box_{-k} u},$  we will make no distinction between $u$ and $\bar{u}$ and write the nonlinearity as $\vec{\lambda} \cdot \nabla u^{2\kappa+1}$. Applying the smoothness of $\mathscr{A}$ and the Strichartz-smoothness estimate,
\begin{align}
&  \|\partial_{x_1} \mathscr{A} ( u^{2\kappa+1} ) \|^{\rm sm}_1 \nonumber\\
&  \lesssim   \sum_{k\in \mathbb{Z}^n, \, |k_1|\ge 20 \vee \max_{j\neq 1} |k_j| } \langle k_1\rangle^{1/2+1/2\kappa} \|\partial_{x_1}\Box_k  \mathscr{A} ( u^{2\kappa+1}  )\|_{L^\infty_{x_1}L^2_{\bar{x} ,t}}  \nonumber\\
 & \lesssim    \sum_{  \mathbb{A}^{2\kappa+1, 1}_{\rm hi} } \  \sum_{k\in \mathbb{Z}^n, \, |k_1|\ge 20 \vee \max_{j\neq 1}|k_j|} \langle k_1\rangle^{1/2+ 1/2\kappa} \left\|\partial_{x_1}\mathscr{A}\Box_k \left( \Box_{k^{(1)}} u... \Box_{k^{(2\kappa+1)}} u \right) \right\|_{L^\infty_{x_1}L^2_{\bar{x} ,t}}  \nonumber \\
  & \quad +\sum_{  \mathbb{A}^{2\kappa+1, 1}_{\rm lo} }  \  \sum_{k\in \mathbb{Z}^n, \, |k_1|\ge 20 \vee \max_{j\neq 1}|k_j|} \langle k_1\rangle^{1/2+ 1/2\kappa} \left\|\partial_{x_1}\mathscr{A}\Box_k \left(\Box_{k^{(1)}} u... \Box_{k^{(2\kappa+1)}} u \right) \right\|_{L^\infty_{x_1}L^2_{\bar{x} ,t}}  \nonumber \\
  & \lesssim    \sum_{  \mathbb{A}^{2\kappa+1, 1}_{\rm hi} } \  \sum_{k\in \mathbb{Z}^n, \, |k_1|\ge 20 \vee \max_{j\neq 1}|k_j|} \langle k_1\rangle^{1/2+ 1/2\kappa} \left\| \Box_k \left( \Box_{k^{(1)}} u... \Box_{k^{(2\kappa+1)}} u \right) \right\|_{L^1_{x_1}L^2_{\bar{x} ,t}}  \nonumber \\
  & \quad +\sum_{  \mathbb{A}^{2\kappa+1, 1}_{\rm lo} }  \  \sum_{k\in \mathbb{Z}^n, \, |k_1|\ge 20 \vee \max_{j\neq 1}|k_j|} \langle k_1\rangle^{1 + 1/2\kappa} \left\| \Box_k \left(\Box_{k^{(1)}} u... \Box_{k^{(2\kappa+1)}} u \right) \right\|_{L^{2(\kappa+1)/(2\kappa+1)}_{x,t}}  \nonumber \\
&:= I+II. \label{ndestI-II}
\end{align}
Using Lemma \ref{boxkest}, we see that in $I$ and $II$, the summation on $k\in \mathbb{Z}^n$ is finitely many and we have the the following restriction on $k\in \mathbb{Z}^n$ in $I$ and $II$:
\begin{align}
|k- k^{(1)}-...- k^{(2\kappa+1)}|_\infty \le 2\kappa + 2. \label{orthog}
\end{align}
So, $|k_1| \le 4(\kappa+1) \max_{s} |k^{(s)}_1|$. We separate the estimate of $I$ into several steps.

{\bf Step 1.} We assume that $|k^{(2)}_1|= \max_{s} |k^{(s)}_1|$.

{\it Case a.}  $|k^{(2)}_1|= \max_{\lambda=1,...,n} |k^{(2)}_\lambda|$. Applying the smooth effect of $\mathscr{A}$ and \eqref{orthog},  we have
\begin{align}
I
 & \lesssim    \sum_{  \mathbb{A}^{2\kappa+1, 1}_{\rm hi} } \ \langle k^{(2)}_1\rangle^{1/2+ 1/2\kappa} \left\|\Box_k \left(\Box_{k^{(1)}} u... \Box_{k^{(2\kappa+1)}} u \right) \right\|_{L^1_{x_1}L^2_{\bar{x} ,t}}  \nonumber \\
& \lesssim     \sum_{k^{(2)}\in \mathbb{Z}^n} \langle k^{(2)}_1 \rangle^{1/2+ 1/2\kappa}  \| \Box_{k^{(2)}} u \|_{L^\infty_{x_1} L^2_{\bar{x} ,t}}   \prod_{s \neq 2}  \sum_{k^{(s)} \in \mathbb{Z}^n }    \| \Box_{k^{(s)}} u \|_{L^{2\kappa}_{x_1}L^\infty_{\bar{x} ,t}}  \nonumber \\
  & \lesssim  \|u\|^{\rm sm} (\|u\|^{\rm max})^{2\kappa}. \label{ndest-1}
\end{align}

{\it Case b.}  $|k^{(2)}_2|= \max_{s=1,...,n} |k^{(2)}_s|$. We have $\langle k_1\rangle \lesssim \langle k^{(2)}_1\rangle  \le  \langle k^{(2)}_2 \rangle$ and
\begin{align}
        \| \Box_{k^{(1)}} u ... \Box_{k^{(2\kappa+1)}} u \|_{L^1_{x_1} L^2_{\bar{x} ,t}}
 & \leqslant  \| \Box_{k^{(1)}} u ... \Box_{k^{(\kappa+1)}} u \|_{ L^2_{x ,t}}
\prod^{2\kappa+1}_{s =\kappa+ 2}    \| \Box_{k^{(s)}} u \|_{L^{2\kappa}_{x_1}L^\infty_{\bar{x} ,t}}  \nonumber \\
  & \lesssim    \| \Box_{k^{(2)}} u \|_{L^{\infty}_{x_2}L^2_{(x_j)_{j\neq 2} ,t}}   \prod^{\kappa+1}_{s=1, s\neq 2}    \| \Box_{k^{(s)}} u \|_{L^{2\kappa}_{x_2}L^\infty_{(x_j)_{j\neq 2} ,t}} \nonumber\\
   & \quad \  \  \times \prod^{2\kappa+1}_{s =\kappa+ 2}    \| \Box_{k^{(s)}} u \|_{L^{2\kappa}_{x_1}L^\infty_{\bar{x} ,t}}. \label{ndest-2}
\end{align}
Hence, noticing that $20 \le |k^{(2)}_1| \le |k^{(2)}_2|$, we have
\begin{align}
I  & \lesssim  \|u\|^{\rm sm} (\|u\|^{\rm max})^{2\kappa}. \label{ndest-3}
\end{align}

{\it Case c.} If $|k^{(2)}_i|= \max_{s=1,...,n} |k^{(2)}_s|$ for some $i>2$, then we can repeat the above proof to get the conclusion.

{\bf Step 2.} We assume that $|k^{(i)}_1|= \max_{s} |k^{(s)}_1|$ for some $i\neq 2$. The estimate of $I$ is the same as in Step 1.
\\

Now we estimate $II$. In view of the Strichartz-smoothness estimate and \eqref{orthog}, we have
\begin{align}
II
 & \lesssim    \sum_{  \mathbb{A}^{2\kappa+1, 1}_{\rm lo} } \ \left\|\Box_{k^{(1)}} u ... \Box_{k^{(2\kappa+1)}} \right\|_{L^{(2\kappa+2)/(2\kappa+1)}_{x,t}} \lesssim (\|u\|^{\rm str})^{2\kappa+1}. \label{ndest-10}
\end{align}
Collecting \eqref{ndest-1}, \eqref{ndest-3} and \eqref{ndest-10}, we have the result, as desired. $\hfill \Box$

\begin{lem}\label{maximalfunctionnorm}
Let $m=2\kappa$ with $\kappa \in \mathbb{N}$.   We have
\begin{align} \label{sm11}
\| \mathcal{T} u\|^{\rm max}  \lesssim  \|u_0\|_{M^{1/2\kappa}_{2,1}} + \| u\|^{\rm sm} (\| u\|^{\rm max})^{2\kappa} + (\| u\|^{\rm str})^{2\kappa+1}
\end{align}
\end{lem}
{\bf Proof.} By symmetry  of $\|\cdot\|^{\rm max}$ it suffices to bound $\| \mathcal{T} u\|^{\rm max}_1$. Applying the maximal function estimate of $S(t)$, one has that
 \begin{align} \label{sm11}
\| \mathcal{T} u\|^{\rm max}_{1}  \lesssim  \|u_0\|_{M^{1/2\kappa}_{2,1}} +  \sum^n_{i=1} \|\partial_{x_i}\mathscr{A} (u^{2\kappa+1} )\|^{\rm max}_1.
\end{align}
We divide the proof into the following two steps.

{\bf Step 1.} $\kappa=1$. We have
\begin{align}
 \|\partial_{x_i} \mathscr{A} ( u^{3} ) \|^{\rm max}_1
& \lesssim    \sum_{  \mathbb{A}^{3, i}_{\rm hi} } \  \sum_{k\in \mathbb{Z}^n}
 \left\|\partial_{x_i}\mathscr{A}\Box_k \left(\Box_{k^{(1)}} u \Box_{k^{(2)}} u\Box_{k^{(3)}} u \right) \right\|_{L^2_{x_1}L^\infty_{\bar{x} ,t}}  \nonumber \\
  & \quad +\sum_{  \mathbb{A}^{3, i}_{\rm lo} }  \  \sum_{k\in \mathbb{Z}^n}  \left\|\partial_{x_i}\mathscr{A}\Box_k \left( \Box_{k^{(1)}} u \Box_{k^{(2)}} u\Box_{k^{(3)}} u\right) \right\|_{L^2_{x_1}L^\infty_{\bar{x} ,t}}  \nonumber \\
  &:= I+II. \label{maxestIII}
\end{align}
In view of Proposition \ref{maxsmest} and \eqref{orthog},
\begin{align}
 I
& \lesssim    \sum_{  \mathbb{A}^{3, i}_{\rm hi} } \  \left\langle \max_{\lambda =1,2,3} |k^{(\lambda)}_i| \right\rangle
 \left\| \Box_{k^{(1)}} u \Box_{k^{(2)}} u\Box_{k^{(3)}} u \right\|_{L^1_{x_1}L^2_{\bar{x} ,t}}.
\end{align}
Assume that $|k^{(s)}_{i}| = |k^{(1)}_{i}| \vee |k^{(2)}_{i}| \vee |k^{(3)}_{i}|.$ By H\"older's inequality,
\begin{align}
 \left\| \Box_{k^{(1)}} u \Box_{k^{(2)}} u\Box_{k^{(3)}} u \right\|_{L^1_{x_1}L^2_{\bar{x} ,t}} \le
 \left\| \Box_{k^{(r)}} u  \right\|_{L^2_{x_1}L^\infty_{\bar{x} ,t}} \left\| \Box_{k^{(s)}} u  \Box_{k^{(q)}} u   \right\|_{L^2_{x ,t}},
\end{align}
where $r,s,q\in \{1,2,3\}$ are different from each other. We can further assume that $|k^{(s)}_{i_0}| = \max_{1\le \lambda \le n}|k^{(s)}_{\lambda}| $. Hence,
\begin{align}
 \left\| \Box_{k^{(1)}} u ... \Box_{k^{(3)}} u \right\|_{L^1_{x_1}L^2_{\bar{x} ,t}} \le
 \left\| \Box_{k^{(r)}} u  \right\|_{L^2_{x_1}L^\infty_{\bar{x} ,t}} \| \Box_{k^{(s)}} u \|_{L^\infty_{x_{i_0}}L^2_{(x_j)_{j\neq i_0} ,t}}
  \| \Box_{k^{(q)}} u  \|_{L^2_{x_{i_0}}L^\infty_{(x_j)_{j\neq i_0} ,t}}.
\end{align}
Noticing that $ 20 \le |k^{(s)}_i | \le | k^{(s)}_{i_0} |$, we have
\begin{align} \label{sm11}
I  \lesssim  \| u\|^{\rm sm} (\| u\|^{\rm max})^{2}.
\end{align}
Using the Strichartz-maximal estimate,
\begin{align}
II
 & \lesssim    \sum_{  \mathbb{A}^{3, i}_{\rm lo} } \ \left\|\Box_{k^{(1)}} u ... \Box_{k^{(3)}} \right\|_{L^{4/3}_{x,t}} \lesssim (\|u\|^{\rm str})^{3}. \label{ndest-100}
\end{align}

{\bf Step 2.} $\kappa\ge 2.$  Using Lemma \ref{connorm},
\begin{align}
 \|\partial_{x_i} \mathscr{A} ( u^{2\kappa+1} ) \|^{\rm max}_1
&  \lesssim   \sum_{k\in \mathbb{Z}^n } \|\partial_{x_i}\Box_k  \mathscr{A} ( u^{2\kappa+1}  )\|_{L^{2\kappa}_{x_1}L^\infty_{\bar{x} ,t}}  \nonumber\\
&  \lesssim
\sum_{k\in \mathbb{Z}^n, |k_1|= \max_{j\neq 1} |k_j| } \|\partial_{x_1}\Box_k  \mathscr{A} ( u^{2\kappa+1}  )\|_{L^{2\kappa}_{x_1}L^\infty_{\bar{x} ,t}}   \nonumber\\
& \quad \quad +...+  \sum_{k\in \mathbb{Z}^n, |k_n|= \max_{j\neq 1} |k_j| } \|\partial_{x_n}\Box_k  \mathscr{A} ( u^{2\kappa+1}  )\|_{L^{2\kappa}_{x_1}L^\infty_{\bar{x} ,t}}   \nonumber\\
& := \Gamma_1+... + \Gamma_n.
\label{maxestIIIhigh}
\end{align}
We estimate $\Gamma_2$ for instance.  By Proposition \ref{str-sm-relation}, we have
\begin{align}
\Gamma_2 &
 \lesssim    \sum_{  \mathbb{A}^{2\kappa+1, 2}_{\rm hi} } \  \sum_{k\in \mathbb{Z}^n, \, |k_2|= \max_{j\neq 1}|k_j|} \langle k_2\rangle^{1/2+ 1/2\kappa} \left\| \Box_k \left(\Box_{k^{(1)}} u ... \Box_{k^{(2\kappa+1)}} u \right) \right\|_{L^1_{x_i}L^2_{(x_j)_{j\neq i} ,t}}  \nonumber \\
  & \quad +\sum_{  \mathbb{A}^{2\kappa+1, 2}_{\rm lo} }  \  \sum_{k\in \mathbb{Z}^n, \, |k_2|= \max_{j\neq 1}|k_j|} \langle k_2\rangle^{1+ 1/2\kappa} \left\|\Box_k \left(\Box_{k^{(1)}} u ... \Box_{k^{(2\kappa+1)}} u \right) \right\|_{L^{2(\kappa+1)/(2\kappa+1)}_{x,t}}  \nonumber \\
  &:= \Gamma_{21} + \Gamma_{22}.
\end{align}
So, using the same way as in the above, we have
\begin{align}
\Gamma_{21} & \lesssim \|u\|^{\rm sm} (\|u\|^{\rm max})^{2\kappa}, \ \  \Gamma_{22} \lesssim (\|u\|^{\rm str})^{2\kappa+1}.
\end{align}
The other terms $\Gamma_i$ can be estimated in an analogous way and the details are omitted. $\quad \Box$

\begin{lem}\label{strichartznorm}
Let $m=2\kappa$ with $\kappa \in \mathbb{N}$.   We have
\begin{align} \label{sm11}
\| \mathcal{T} u\|^{\rm str}  \lesssim  \|u_0\|_{M^{1/2\kappa}_{2,1}} + \| u\|^{\rm sm} (\| u\|^{\rm max})^{2\kappa} + (\| u\|^{\rm str})^{2\kappa+1}
\end{align}
\end{lem}
{\bf Proof.} By Strichartz estimate, we have
\begin{align} \label{sm11}
\| \mathcal{T} u\|^{\rm str}  \lesssim  \|u_0\|_{M^{1/2\kappa}_{2,1}} + \sum^n_{i=1} \|\partial_{x_i} \mathscr{A}
u^{2 \kappa +1}\|^{\rm str}.
\end{align}
It suffices to bound  $ \|\partial_{x_1} \mathscr{A}
u^{2 \kappa +1}\|^{\rm str}$. In view of the Strichartz and smooth-Strichartz estimates,
\begin{align}
\|\partial_{x_1} \mathscr{A}
u^{2 \kappa +1}\|^{\rm str}
&\lesssim  \sum_{k \in \mathbb{Z}^n } \langle k_1\rangle^{1/2+ 1/2\kappa}
\sum_{\mathbb{A}^{2\kappa+1,1}_{\rm hi}} \|\Box_k \left(\Box_{k^{(1)}} u ...
\Box_{k^{(2\kappa+1)}} u \right) \|_{L^{1}_{x_1}
L^2_{\bar{x},t}}  \nonumber\\
& \ \  \  +  \sum_{k \in \mathbb{Z}^n} \langle
k_1\rangle^{1+1/2\kappa} \sum_{\mathbb{A}^{2\kappa+1,1}_{\rm lo}} \| \Box_k
\left(\Box_{k^{(1)}} u ... \Box_{k^{(2\kappa +1)}} u
\right)\|_{L^{2(1+\kappa)/(1+2\kappa)}_{t,x}  }.
\label{strestf}
\end{align}
Repeating the argument as in Lemma \ref{smthnorm}, we obtain the result and the details are omitted. $\hfill \Box$

\section{ Cubic nonlinearity in 2D }

Now we split the (semi-)norms in different directions.
\begin{align}
 & \|u\|^{\rm sm2}_1= \sum_{k\in \mathbb{Z}^2, \, |k_1|\ge |k_2|\vee 20} \langle k_1\rangle^{3/2} \|\Box_k u\|_{L^\infty_{x_1}L^2_{x_2,t}},  \\
 & \|u\|^{\rm sm2}_2= \sum_{k\in \mathbb{Z}^2, \, |k_2|\ge |k_1|\vee 20} \langle k_2\rangle^{3/2} \|\Box_k u\|_{L^\infty_{x_2}L^2_{x_1,t}};
\end{align}
\begin{align}
  \|u\|^{\rm max}_1= \sum_{k\in \mathbb{Z}^2}   \|\Box_k u\|_{L^2_{x_1}L^\infty_{x_2,t}}, \ \  \|u\|^{\rm max}_2= \sum_{k\in \mathbb{Z}^2}   \|\Box_k u\|_{L^2_{x_2}L^\infty_{x_1,t}};
\end{align}
\begin{align}
  \|u\|^{\rm ant}_1= \sum_{k\in \mathbb{Z}^2}   \|\Box_k u\|_{L^2_{x_1}L^4_{x_2,t}}, \ \  \|u\|^{\rm ant}_2= \sum_{k\in \mathbb{Z}^2}   \|\Box_k u\|_{L^2_{x_2}L^4_{x_1,t}};
\end{align}
We see that
$$
\|u\|^{\rm sm2}=\|u\|^{\rm sm2}_1+\|u\|^{\rm sm2}_2, \ \  \|u\|^{\rm max}=\|u\|^{\rm max}_1+\|u\|^{\rm max}_2, \ \ \|u\|^{\rm ant}=\|u\|^{\rm ant}_1+\|u\|^{\rm ant}_2.
$$
Due to $\Box_k \bar{u} = (-1)^n \overline{\Box_{-k} u}$, we can assume that
$$
F(u, \bar{u}, \nabla u, \nabla \bar{u})= F(u,  \nabla u)= \sum^\infty_{j=3} \ \sum_{\kappa+\nu_1+\nu_2=j} C_{\kappa\nu_1 \nu_2} u^\kappa u^{\nu_1}_{x_1} u^{\nu_2}_{x_2}:  = \sum^\infty_{j=3} F_j (u,  \nabla u).
$$
We define the following
\begin{align}
\mathcal{D} = \left\{u: \sum_{\alpha=0,1} \sum_{i=1,2}\|\partial^\alpha_{x_i}u\|^{\rm sm2 \,\cap\, max \,\cap\, ant \,\cap\, str2 \,\cap\, gstr}  \le \delta \right\},
\end{align}
and for any $u,v\in \mathcal{D}$,
\begin{align}
 d(u,v)= \sum_{\alpha=0,1} \sum_{i=1,2} \|\partial^\alpha_{x_i}(u-v)\|^{\rm sm2 \,\cap\, max \,\cap\, ant \,\cap\, str2 \,\cap\, gstr}.
\end{align}
We consider the following mapping $\mathscr{T}$ in $(\mathcal{D}, d)$,
 \begin{align}
\mathscr{T}: u\to S(t) u_0 - {\rm i}  \sum^\infty_{j=3}  \int^t_0 S(t-\tau) F_j (u(\tau), \nabla u(\tau))d\tau,
\end{align}
 By Lemma \ref{connorm},
\begin{lem}\label{simpleseminorm}
We have
\begin{align} \label{ssn}
 \sum_{i=1,2}\sum_{\alpha=0,1}\|\partial^\alpha_{x_i}u\|^{\rm sm2} \lesssim \|\partial_{x_1}u\|^{\rm sm2}_1 + \|\partial_{x_2}u\|^{\rm sm2}_2.
\end{align}
\end{lem}
Using Lemma \ref{simpleseminorm}, we have
\begin{align}
 \sum_{\alpha=0,1} \sum_{i=1,2}\|\partial^\alpha_{x_i} \mathscr{T}u\|^{\rm sm2} \lesssim \|\partial_{x_1}\mathscr{T} u\|^{\rm sm2}_1 + \|\partial_{x_2}\mathscr{T} u\|^{\rm sm2}_2. \nonumber
\end{align}

\begin{lem} \label{smest}
Let $u\in \mathcal{D}$. Then we have
\begin{align}
   \|\partial_{x_1}\mathscr{T} u\|^{\rm sm2}_1 + \|\partial_{x_2}\mathscr{T} u\|^{\rm sm2}_2 \lesssim \|u_0\|_{M^2_{2,1}} +    \sum^\infty_{j=3} C^j \delta^{ j} . \nonumber
\end{align}
\end{lem}
{\bf Proof.} First, we estimate $ \|\partial_{x_1}\mathscr{T} u\|^{\rm sm2}_1$.
\begin{align}
 \|\partial_{x_1}\mathscr{T} u\|^{\rm sm2}_1  \le \|\partial_{x_1} S(t) u_0\|^{\rm sm2}_1 +  \sum^\infty_{j=3} \|\partial_{x_1} \mathscr{A} F_j(u, \nabla u)\|^{\rm sm2}_1.  \nonumber
\end{align}
By Proposition \ref{smeff},
\begin{align}
 \|\partial_{x_1} S(t) u_0\|^{\rm sm2}_1
 & \lesssim  \sum_{k\in \mathbb{Z}^2, \, |k_1|\ge |k_2|\vee 20} \langle k_1\rangle^2 \| \Box_k u_0\|_{2} \le \|u_0\|_{M^2_{2,1}}.  \nonumber
\end{align}
For convenience, we write
\begin{align}
& \mathbb{A}^{\lambda, i}_{\rm lo} = \left\{(k^{(1)},..., k^{(\lambda)}) \in (\mathbb{Z}^2)^\lambda: \, \max_{1\le s\le \lambda} |k^{(s)}_i| < 20 \right\}, \nonumber\\
& \mathbb{A}^{\lambda, i}_{\rm hi} = \left\{(k^{(1)},..., k^{(\lambda)}) \in (\mathbb{Z}^2)^\lambda: \, \max_{1\le s\le \lambda} |k^{(s)}_i| \ge 20 \right\},
\end{align}
where $k^{(s)}=(k^{(s)}_1, k^{(s)}_2)$. Let us write $v_1=...=v_\kappa =u$, $v_{\kappa+1}=...=v_{\kappa+ \nu_1}= u_{x_1}$ and $v_{\kappa+\nu_1+1}=...=v_{\kappa+ \nu_1+ \nu_2}= u_{x_2}$,
$$
F_j(u, \nabla u)= v_1...v_{\kappa+|\nu|}, \ \ |\nu| = \nu_1+ \nu_2.
$$
We see that for $\kappa+|\nu|=j$,
\begin{align}
&  \|\partial_{x_1} \mathscr{A} (v_1...v_{\kappa+|\nu|})  \|^{\rm sm2}_1 \nonumber\\
 & \lesssim   \left(\sum_{  \mathbb{A}^{j, 1}_{\rm lo} }+ \sum_{  \mathbb{A}^{j, 1}_{\rm hi} } \right) \sum_{k\in \mathbb{Z}^2, \, |k_1|\ge |k_2|\vee 20} \langle k_1\rangle^{3/2}
    \left\|\partial_{x_1}\mathscr{A}\Box_k \left(\prod^j_{s=1}\Box_{k^{(s)}} v_s    \right) \right\|_{L^\infty_{x_1}L^2_{x_2,t}} \nonumber\\
  &:= I+II. \label{estI-II}
\end{align}
By Proposition \ref{str-sm-relation}, H\"older's inequality and Lemma \ref{boxkest},
\begin{align}
 I
 & \lesssim    \sum_{ \mathbb{A}^{j, 1}_{\rm lo} }  \    \sum_{k\in \mathbb{Z}^2, \, |k_1|\ge |k_2|\vee 20} \langle k_1\rangle^2
   \left\| \Box_k \left(\prod^j_{s=1}\Box_{k^{(s)}} v_s   \right) \right\|_{L^{4/3}_{x,t}} \nonumber\\
& \lesssim  \sum_{ \mathbb{A}^{j, 1}_{\rm lo} }  \  \sum_{k\in \mathbb{Z}^2, \, |k_1|\ge |k_2|\vee 20} \langle k_1\rangle^2  \prod^{j-3}_{s=1} \|\Box_{k^{(s)}} v_s \|_{L^\infty_{x,t}} \prod^{j }_{s=j-2} \|\Box_{k^{(s)}} v_s \|_{L^4_{x,t}} \nonumber \\
 &  \quad\quad  \times \prod_{i=1,2}\chi_{\left(|k_i-k^{(1)}_i-...-k^{(j)}_i|\le  j+1 \right)}. \label{estI1}
\end{align}
Noticing that there are at most $O(j)$ non-zero terms in the summation $\sum_{k\in \mathbb{Z}^2, \, |k_1|\ge |k_2|\vee 20}$ and $|k_1| \le Cj$ in \eqref{estI1}, we easily see that
\begin{align}
 I
 & \lesssim    C^j  \sum_{ \mathbb{A}^{j, 1}_{\rm lo} } \|\Box_{k^{(s)}} v_s \|_{L^\infty_{x,t}} \prod^{j }_{s=j-2} \|\Box_{k^{(s)}} v_s \|_{L^4_{x,t}} \lesssim     C^j  (\| u\|^{\rm str2})^{j} . \label{estI2}
\end{align}
For convenience, we further write
\begin{align}
& \mathbb{A}^{\lambda, i}_{\ell, \, \rm hi} = \left\{(k^{(1)},..., k^{(\lambda)}) \in \mathbb{A}^{\lambda, i}_{ \rm hi}: \, k^{(\ell )}_{i} = \max_{1\le s\le \lambda} |k^{(s)}_i|   \right\}, \ \  \ell =1,...,\lambda.
\end{align}
By Proposition \ref{smeff},
\begin{align}
 II
 & \lesssim    \left(\sum_{ \mathbb{A}^{j, 1}_{1, \rm hi} }   +...+ \sum_{ \mathbb{A}^{j, 1}_{j, \rm hi} } \right) \sum_{k\in \mathbb{Z}^2, \, |k_1|\ge |k_2|\vee 20} \langle k_1\rangle^{3/2}
    \left\| \Box_k \left(\prod^j_{s=1}\Box_{k^{(s)}} v_s    \right) \right\|_{L^1_{x_1}L^2_{x_2,t}} \nonumber\\
&  :=  II_1 +...+ II_{j} .   \label{estII1}
\end{align}
The estimates for $ II_1,...,  II_{j} $ are similar and it suffices to estimate $ II_{j-2}$.  If $|k^{(j-2)}_{1}| = \max_{1\le s\le j} |k^{(s)}_1| $,  by H\"older's inequality and Lemma \ref{boxkest}, we have
\begin{align}
& \left\| \Box_k \left( \Box_{k^{(1)}} v_1 ...  \Box_{k^{(j)}} v_j   \right) \right\|_{L^1_{x_1}L^2_{x_2,t}} \nonumber\\
&   \le \prod^{j-3}_{s=1}   \|\Box_{k^{(s)}} v_s \|_{L^\infty_{x,t}} \|\Box_{k^{(j-1)}} v_{j-1}  \|_{L^2_{x_1}L^\infty_{x_2,t}} \|\Box_{k^{(j-2)}} v_{j-2} \Box_{k^{(j)}} v_{j}\|_{L^2_{x,t}}  \nonumber\\
&  \quad\quad \times \prod_{i=1,2}\chi_{\left(|k_i-k^{(1)}_i-...-k^{(j)}_i|\le  j+1 \right)}.   \label{estII2}
\end{align}
In \eqref{estII2}, if $|k^{(j-2)}_1| \ge |k^{(j-2)}_2|$, by H\"older's inequality,
\begin{align}
  \|\Box_{k^{(j-2)}} v_{j-2} \Box_{k^{(j)}} v_{j}\|_{L^2_{x,t}}  \le  \|\Box_{k^{(j-2)}} v_{j-2} \|_{L^\infty_{x_1}L^2_{x_2,t}}  \| \Box_{k^{(j)}} v_{j}\|_{L^2_{x_1}L^\infty_{x_2,t}};   \label{estII3}
\end{align}
and  if $|k^{(j-2)}_2| > |k^{(j-2)}_1|$, by H\"older's inequality,
\begin{align}
  \|\Box_{k^{(j-2)}} v_{j-2} \Box_{k^{(j)}}v_{j }\|_{L^2_{x,t}}  \le  \|\Box_{k^{(j-2)}} v_{j-2} \|_{L^\infty_{x_2}L^2_{x_1,t}}  \| \Box_{k^{(j)}} v_{j }\|_{L^2_{x_2}L^\infty_{x_1,t}}.   \label{estII4}
\end{align}
We see that in $II_{j-2}$, there are at most $Cj$ non-zeor terms in the summation $\sum_{k\in \mathbb{Z}^2, \, |k_1|\ge |k_2|\vee 20}$  and $|k_1| \le Cj (|k^{(j-2)}_1|\vee |k^{(j-2)}_2|)$ if $k^{(1)},...,k^{(j)} \in \mathbb{A}^{j, 1}_{\rm hi}$, it follows that
\begin{align}
 II_{j-2}  \lesssim   C^j  \left(\prod^{j-3}_{s=1} \|v_s\|^{\rm str2} \right) \|v_{j-2}\|^{\rm sm2}   \|v_{j-1}\|^{\rm max}  \|v_{j}\|^{\rm max}  \label{estII5}
\end{align}
So, we have
\begin{align}
 \|\partial_{x_1} \mathscr{A} (v_1...v_{\kappa+|\nu|})  \|^{\rm sm2}_1   \lesssim      C^j \delta^{ j}   \label{estII10}
\end{align}
holds for any $u\in \mathcal{D}$. $\hfill\Box$

\begin{lem} \label{maxantest}
Let $u\in \mathcal{D}$. Then we have
\begin{align}
  \sum_{\alpha=0,1} \sum_{i=1,2} \|\partial^\alpha_{x_i}\mathscr{T} u\|^{\rm max \, \cap \, ant} \lesssim \|u_0\|_{M^2_{1,1}} +    \sum^\infty_{j=3} C^j \delta^{ j} . \nonumber
\end{align}
\end{lem}
{\bf Proof.} We have
\begin{align}
 & \sum_{\alpha=0,1} \sum_{i=1,2} \|\partial^\alpha_{x_i}\mathscr{T} u\|^{\rm max \, \cap \, ant} \nonumber\\
  &  \le \sum_{\alpha=0,1} \sum_{i=1,2} \|\partial^\alpha_{x_i} S(t) u_0\|^{\rm max \, \cap \, ant} +  \sum^\infty_{j= 3} \sum_{\alpha=0,1} \sum_{i=1,2} \|\partial^\alpha_{x_i}  \mathscr{A} F_j (u, \nabla u)\|^{\rm max \, \cap \, ant} \nonumber
\end{align}
By Corollary \ref{loc-gfexpress}, we have for $\alpha=0,1$ and $i=1,2$,
\begin{align}
 \|\partial^\alpha_{x_i} S(t) u_0\|^{\rm max \, \cap \, ant}_1  & \leq \sum_{k\in \mathbb{Z}^2}  \max_{p=\infty,\, 4} \|\partial^\alpha_{x_i}\Box_k S(t) u_0\|_{L^2_{x_1}L^p_{x_2,t}}  \nonumber\\
 & \lesssim  \sum_{k\in \mathbb{Z}^2} \langle k_1\rangle^{1/2} \langle k_i\rangle \| \Box_k u_0\|_{1} \le \|u_0\|_{M^2_{1,1}}.  \nonumber
\end{align}
Again, in view of Corollary \ref{loc-gfexpress},
\begin{align}
& \sum_{\alpha=0,1} \sum_{i=1,2} \|\partial^\alpha_{x_i}  \mathscr{A} (v_1... v_j) \|_1^{\rm max \, \cap \, ant} \nonumber\\
& \le   \sum_{\alpha=0,1} \sum_{i=1,2} \sum_{k\in \mathbb{Z}^2} \|\partial^\alpha_{x_i}\Box_k  \mathscr{A} ( v_1... v_j )\|_{L^2_{x_1}L^4_{x_2,t} \cap L^2_{x_1}L^\infty_{x_2,t}}  \nonumber\\
& \le    \sum_{k\in \mathbb{Z}^2, |k_1|\ge |k_2|}  \ \sum_{k^{(1)},...,k^{(j)}\in \mathbb{Z}^{2}} \langle k_1\rangle^{3/2}
  \left\|\Box_k \left( \Box_{k^{(1)}} v_1 ... \Box_{k^{(j)}} v_j    \right) \right\|_{L^1_{x,t}} \nonumber\\
& \quad  \quad  +     \sum_{k\in \mathbb{Z}^2, |k_2|\ge |k_1|}  \ \sum_{k^{(1)},...,k^{(j)}\in \mathbb{Z}^{2}} \langle k_2\rangle^{3/2}
  \left\|\Box_k \left(\Box_{k^{(1)}} v_1 ... \Box_{k^{(j)}} v_j \right) \right\|_{L^1_{x,t}} \nonumber\\
& := \Upsilon_1 + \Upsilon_2.  \label{upsilon}
\end{align}
We now estimate $\Upsilon_1$. We have
\begin{align}
 \Upsilon_1 & \le    \sum_{k\in \mathbb{Z}^2, |k_1|\ge |k_2|}  \ \sum_{ \mathbb{A}^{j,1}_{\rm lo}} \langle k_1\rangle^{3/2}
  \left\|\Box_k \left(\Box_{k^{(1)}} v_1 ... \Box_{k^{(j)}} v_j \right) \right\|_{L^1_{x,t}} \nonumber\\
&  \quad +    \sum_{k\in \mathbb{Z}^2, |k_1|\ge |k_2|}  \ \sum_{ \mathbb{A}^{j,1}_{\rm hi}} \langle k_1\rangle^{3/2}
  \left\|\Box_k \left(\Box_{k^{(1)}} v_1 ... \Box_{k^{(j)}} v_j \right) \right\|_{L^1_{x,t}} \nonumber\\
& := \Upsilon_{11} + \Upsilon_{12}.
\end{align}
By H\"older's inequality and Lemma \ref{boxkest},
\begin{align}
 \Upsilon_{11} & \le    C^j     \ \sum_{k^{(1)},...,k^{(j)}\in \mathbb{Z}^2}
  \prod^{j-3}_{s=1}\|\Box_{k^{(s)}} v_s \|_{L^\infty_{x,t}}  \prod^{j}_{s=j-2} \|\Box_{k^{(s)}} v_s \|_{L^3_{x,t}} \nonumber\\
 & \lesssim    C^j \left( \prod^{j-3}_{s=1}
    \| v_s\|^{\rm str2}  \right)   \left(\prod^{j}_{s=j-2} \| v_s\|^{\rm gstr}\right) .
\end{align}
In order to bound $\Upsilon_{12}$, we further decompose $\mathbb{A}^{j,1}_{\rm hi}$. We have
\begin{align}
\Upsilon_{12}
  \le     \sum^{j}_{\ell=1}  \sum_{k\in \mathbb{Z}^2, |k_1|\ge |k_2|}  \ \sum_{\mathbb{A}^{j,1}_{\ell,  \rm hi}} \langle k_1\rangle^{3/2} \left\|\Box_k \left(\prod^j_{s=1}\Box_{k^{(s)}} v_s    \right) \right\|_{L^1_{x,t}}
& :=    \sum^{j}_{\ell=1} \Upsilon_{12, \ell}.
\end{align}
The estimates of $\Upsilon_{12, \ell}$ for $\ell=1,...,j$ are similar and we only need to estimate $\Upsilon_{12, 1}$. By H\"older's inequality and Lemma \ref{boxkest},
\begin{align}
   \Upsilon_{12, 1}
&  \lesssim    C^j   \sum_{ \mathbb{A}^{j,1}_{1, \rm hi}} \langle k^{(1)}_1\rangle^{3/2}
    \times \prod^{j-2}_{s=2}\|\Box_{k^{(s)}} v_s \|_{L^\infty_{x,t}}  \|\Box_{k^{(1)}} v_1  \Box_{k^{(j-1)}} v_{j-1}  \Box_{k^{(j)}}  v_{j}   \|_{L^1_{x,t}} \nonumber\\
&  \lesssim    C^j   \sum_{ \mathbb{A}^{j,1}_{1, \rm hi}, \, |k^{(1)}_1| \ge |k^{(1)}_2|} \langle k^{(1)}_1\rangle^{3/2}  \prod^{j-2}_{s=2}\|\Box_{k^{(s)}} v_s\|_{L^\infty_{x,t}}
  \|\Box_{k^{(1)}} v_1 \|_{L^\infty_{x_1}L^2_{x_2,t}}   \nonumber\\
   & \quad\quad \times \|\Box_{k^{(j-1)}}  v_{j-1} \|_{L^2_{x_1}L^4_{x_2,t}}  \|\Box_{k^{(j)}}  v_{j}   \|_{L^2_{x_1}L^4_{x_2,t}} \nonumber\\
 & \quad +    C^j   \sum_{ \mathbb{A}^{j,1}_{1, \rm hi}, \, |k^{(1)}_1| < |k^{(1)}_2|} \langle k^{(1)}_2\rangle^{3/2}  \prod^{j-2}_{s=2}\|\Box_{k^{(s)}} v_s \|_{L^\infty_{x,t}}
 \|\Box_{k^{(1)}} v_1  \|_{L^\infty_{x_2}L^2_{x_1,t}}   \nonumber\\
 & \quad\quad  \times  \|\Box_{k^{(j-1)}}  v_{j-1} \|_{L^2_{x_2}L^4_{x_1,t}}  \|\Box_{k^{(j)}}  v_{j}  \|_{L^2_{x_2}L^4_{x_1,t}} \nonumber\\
 & \lesssim    C^j \|v_1\|^{\rm sm2} \left( \prod^{j-2}_{s=2}
    \|v_s\|^{\rm str2} \right)  \prod^{j}_{s=j-1} \| v_s \|^{\rm ant} .
 \end{align}
Analogous to $\Upsilon_{12,1}$,  we have
\begin{align}
 \Upsilon_{12, \ell}    \lesssim     C^j  \delta^{ j}.
 \end{align}
Using the same way as in the estimates of  $\Upsilon_{1}$, we can obtain that
\begin{align}
 \Upsilon_{2}   \lesssim      C^j  \delta^{ j}.
 \end{align}
The result follows. $\hfill\Box$

\begin{lem} \label{gstrest}
Let $u\in \mathcal{D}$. Then we have
\begin{align}
  \sum_{\alpha=0,1} \sum_{i=1,2} \|\partial^\alpha_{x_i}\mathscr{T} u\|^{\rm gstr} \lesssim \|u_0\|_{M^2_{1,1}} +    \sum^\infty_{j=3} C^j \delta^{ j} . \nonumber
\end{align}
\end{lem}
{\bf Proof.} We have
\begin{align}
 \sum_{\alpha=0,1} \sum_{i=1,2} \|\partial^\alpha_{x_i}\mathscr{T} u\|^{\rm gstr }
  &  \le \sum_{\alpha=0,1} \sum_{i=1,2} \|\partial^\alpha_{x_i} S(t) u_0\|^{\rm gstr} + \sum^\infty_{j=3} \sum_{\alpha=0,1} \sum_{i=1,2} \|\partial^\alpha_{x_i}  \mathscr{A} F_j(u, \nabla u)\|^{\rm gstr}. \nonumber
   \end{align}
By Corollary \ref{loc-gfexpress}, we have for $\alpha=0,1$ and $i=1,2$,
\begin{align}
 \|\partial^\alpha_{x_i} S(t) u_0\|^{\rm gstr}
 \lesssim  \sum_{k\in \mathbb{Z}^2} \langle k_1\rangle^{1/3} \langle k_i\rangle \| \Box_k u_0\|_{1} \le \|u_0\|_{M^2_{1,1}}.  \nonumber
\end{align}
Again, in view of Corollary \ref{loc-gfexpress},
\begin{align}
  \sum_{\alpha=0,1} \sum_{i=1,2} \|\partial^\alpha_{x_i}  \mathscr{A}  ( v_1...v_{j} )\|^{\rm gstr}
 & \le   \sum_{k\in \mathbb{Z}^2, \, |k_1|\ge |k_2|} \langle k_1\rangle^{4/3} \| \Box_k    (v_1...v_{j})\|_{L^1_{x,t} }  \nonumber\\
 & \quad +     \sum_{k\in \mathbb{Z}^2, \, |k_2|> |k_1|} \langle k_2\rangle^{4/3} \| \Box_k    (v_1...v_{j})\|_{L^1_{x,t} }  \nonumber\\
& \le \Upsilon_1 + \Upsilon_2,
\end{align}
where $\Upsilon_1$ and $\Upsilon_2$ are the same as in \eqref{upsilon}. So,one can repeat the proof of Lemma \ref{maxantest} to obtain the result. $\hfill\Box$

\begin{lem} \label{strest}
Let $u\in \mathcal{D}$. Then we have
\begin{align}
  \sum_{\alpha=0,1} \sum_{i=1,2} \|\partial^\alpha_{x_i}\mathscr{T} u\|^{\rm str2} \lesssim \|u_0\|_{M^2_{1,1}} +   \sum^\infty_{j=3} C^j \delta^{ j} . \nonumber
\end{align}
\end{lem}
{\bf Proof.} We have
\begin{align}
   \sum_{\alpha=0,1} \sum_{i=1,2} \|\partial^\alpha_{x_i}\mathscr{T} u\|^{\rm str }
  &  \le \sum_{\alpha=0,1} \sum_{i=1,2} \|\partial^\alpha_{x_i} S(t) u_0\|^{\rm str2} +  \sum^\infty_{ j=3} \sum_{\alpha=0,1} \sum_{i=1,2} \|\partial^\alpha_{x_i}  \mathscr{A} F_j(u, \nabla u)\|^{\rm  str}. \nonumber
\end{align}
By Corollary \ref{loc-gfexpress}, we have for $\alpha=0,1$ and $i=1,2$,
\begin{align}
 \|\partial^\alpha_{x_i} S(t) u_0\|^{\rm  str2}
 & \lesssim  \sum_{k\in \mathbb{Z}^2} \langle k \rangle  \langle k_i\rangle \| \Box_k u_0\|_{2} \le \|u_0\|_{M^2_{2,1}} \le  \|u_0\|_{M^2_{1,1}}.  \nonumber
\end{align}
In view of  Proposition \ref{strichartz}, Lemma \ref{boxkest} and  H\"older's inequality,
\begin{align}
 \|  \mathscr{A}  (v_1...v_{j})\|^{\rm str2}
 & \le    \sum_{k\in \mathbb{Z}^2} \langle k \rangle  \| \Box_k    (v_1...v_{j})\|_{ L^{4/3}_{x,t} }  \nonumber\\
& \lesssim   C^j \sum_{k^{(1)},...,k^{(j)}\in \mathbb{Z}^2}    \sum_{k\in \mathbb{Z}^2} \langle |k^{(1)}|\vee...\vee|k^{(j)|} \rangle
    \left\| \Box_k \left(\prod^j_{s=1}\Box_{k^{(s)}} v_s \right) \right\|_{L^{4/3}_{x,t}} \nonumber\\
& \lesssim  C^j  \sum_{k^{(1)},...,k^{(j)}\in \mathbb{Z}^2}    \langle |k^{(1)}|\vee...\vee|k^{(j)|} \rangle  \prod^{ j}_{s=1} \|\Box_{k^{(s)}}v_s \|_{L^\infty_{x,t} \cap L^4_{x,t}} \nonumber \\
& \lesssim    C^j  \prod^j_{s=1} \| v_s\|^{\rm str2} . \label{str1}
\end{align}
Next, we estimate
\begin{align}
  \sum_{i=1,2} \| \partial_{x_i} \mathscr{A}  (v_1...v_{j}) \|^{\rm str2}
 & \lesssim    \sum_{ |k|\le 40}      \| \Box_k    \mathscr{A} (v_1...v_{j} )\|_{L^\infty_tL^2_x \, \cap \,L^4_{x,t} }  \nonumber\\
 & +    \sum_{ |k|> 40} \langle k \rangle \sum_{i=1,2} \| \Box_k    \partial_{x_i} \mathscr{A} (v_1...v_{j})\|_{L^\infty_tL^2_x \, \cap \,L^4_{x,t} }  \nonumber\\
 &:= \Gamma_1+ \Gamma_2.
  \label{str2}
\end{align}
 Using the same way as in \eqref{str1}, $\Gamma_1$ can be estimated in an analogous way as above.
Now we consider the estimate of $\Gamma_2$.
It follows from Lemma \ref{connorm} that
\begin{align}
\Gamma_2
&  \lesssim  \sum^\infty_{j=0}   \sum_{k\in \mathbb{Z}^2,\,  |k_1| \ge |k_2|\vee 20} \langle k_1 \rangle   \| \Box_k    \partial_{x_1} \mathscr{A} (v_1...v_{j})\|_{L^\infty_tL^2_x \, \cap \,L^4_{x,t} }  \nonumber\\
& \quad +   \sum^\infty_{j=0}   \sum_{k\in \mathbb{Z}^2,\,  |k_2| > |k_1|\vee 20} \langle k_2 \rangle   \| \Box_k    \partial_{x_2} \mathscr{A} (v_1...v_{j})\|_{L^\infty_tL^2_x \, \cap \,L^4_{x,t} } \nonumber\\
& := \Gamma_{21} + \Gamma_{22}.
  \label{str5}
\end{align}
We have
\begin{align}
\Gamma_{22}
 & \lesssim   \left(\sum_{ \mathbb{A}^{j, 2}_{\rm lo} }+ \sum_{ \mathbb{A}^{j, 2}_{\rm hi} } \right) \sum_{k\in \mathbb{Z}^2, \, |k_2|\ge |k_1|\vee 20} \langle k_2\rangle
    \left\|\partial_{x_2}\mathscr{A}\Box_k \left(\prod^j_{s=1}\Box_{k^{(s)}} v_s   \right) \right\|_{L^\infty_tL^2_x \cap L^4_{x,t}} \nonumber\\
  &:= \Gamma_{22,1}+ \Gamma_{22,2}.
\end{align}
By Proposition \ref{strichartz},
\begin{align}
  \Gamma_{22,1}
 & \lesssim   \sum_{ \mathbb{A}^{j, 2}_{\rm lo} }    \sum_{k\in \mathbb{Z}^2, \, |k_2|\ge |k_1|\vee 20} \langle k_2\rangle^2
   \left\| \Box_k \left(\prod^j_{s=1}\Box_{k^{(s)}} v_s   \right) \right\|_{L^{4/3}_{x,t}}.
\end{align}
Then we can use the same way as in \eqref{estI1}--\eqref{estI2} to obtain that
\begin{align}
\Gamma_{22,1}
 & \lesssim     C^j  \prod^j_{s=1} \|v_s\|^{\rm str2} . \label{str6}
\end{align}
By Proposition \ref{str-sm-relation},
\begin{align}
 \Gamma_{22,2}
 & \lesssim     \sum_{ \mathbb{A}^{j, 2}_{\rm hi} }    \sum_{k\in \mathbb{Z}^2, \, |k_2|> |k_1|\vee 20} \langle k_2\rangle^{3/2}
   \left\| \Box_k \left(\prod^j_{s=1}\Box_{k^{(s)}} v_s   \right) \right\|_{L^1_{x_2}L^2_{x_1,t}},   \label{str7}
\end{align}
which reduces to the estimate of $II$  as in \eqref{estI-II} and we have
\begin{align}
 \Gamma_{22,2}
 & \lesssim    C^j \delta^{j}.  \label{str8}
\end{align}
Analogous to $\Gamma_{22,2}$,  $\Gamma_{22,1}$ can be bounded by the right hand side of \eqref{str8}.  $\hfill\Box$

\medskip
\medskip

{\bf Proof of Theorem \ref{thg1}.}  By Lemmas \ref{smest}, \ref{maxantest}, \ref{gstrest} and \ref{strest}, we immediately have for any $u\in \mathcal{D}$,
  \begin{align}
  \sum_{\alpha=0,1}\sum_{i=1,2} \|\partial^\alpha_{x_i} \mathscr{T} u\|^{\rm sm2 \, \cap \, max \, \cap \, ant \cap \, str2 \cap \,  gstr} \lesssim \|u_0\|_{M^2_{1,1}} + \delta^3.
\end{align}
Similarly, for any $u,v \in \mathcal{D}$,
 \begin{align}
  d(u,v) \lesssim   \delta^2 d(u,v).
\end{align}
Following a standard contraction mapping argument, we obtain that \eqref{Schmap} has a unique solution $u \in C(\mathbb{R}, M^2_{2,1})  \cap X$.

Finally, it suffices to show that $u \in C_{\rm loc}(\mathbb{R}, M^{3/2}_{1,1})$. Let $T>0$ be arbitrary.  By Proposition \ref{ubinL1},
 \begin{align}
\|u\|_{C([-T,T]; M^{3/2}_{1,1})} & \lesssim \langle T\rangle \|u_0\|_{M^{3/2}_{1,1}} +  \langle T\rangle \int^T_{-T} \|F(u(\tau))\|_{M^{3/2}_{1,1}} d\tau \nonumber\\
& \lesssim \langle T\rangle \|u_0\|_{M^{3/2}_{1,1}} +  \langle T\rangle \sum_{k\in \mathbb{Z}^2} \langle k\rangle^{3/2} \|\Box_k F(u(\tau))\|_{L^1_{x,t}},
\end{align}
which reduces to the estimate as in  Lemma \ref{maxantest} if we treat $T>0$ as a fixed number.
$\hfill \Box$

\medskip

\section{Quartic nonlinearity in 1D}

By Corollary \ref{loc-gfexpress},
 \begin{align}
\|\Box_k S(t) u_0\|_{L_{x}^{3}L^{\infty}_{t} (\mathbb{R}^{1+1} )} & \lesssim \langle k\rangle^{1/3} \|\Box_k u_0\|_{L^1(\mathbb{R})},  \label{fre-local-1d1}\\
\|\Box_k \mathscr{A} f\|_{L_{x }^{3}L^{\infty}_{t} (\mathbb{R}^{1+1} )} & \lesssim \langle k\rangle^{1/3} \|\Box_k f\|_{L^1_{x,t} (\mathbb{R}^{1+1})}. \label{fre-local-1d2}
\end{align}
\begin{align}
 \|\Box_k S(t) u_0\|_{L_{x }^{3}L^{6}_{ t} (\mathbb{R}^{1+1} )} & \lesssim \langle k\rangle^{1/3} \|\Box_k u_0\|_{L^1(\mathbb{R} )},  \label{fre-local-1d3}\\
 \|\Box_k \mathscr{A} f\|_{L_{x }^{3}L^6_{ t} (\mathbb{R}^{1+1} )} & \lesssim \langle k\rangle^{1/3} \|\Box_k f\|_{L^1_{x,t} (\mathbb{R}^{1+1})}. \label{fre-local-1d4}
\end{align}
\begin{align}
\|\Box_k S(t) u_0\|_{L^{4}_{x ,t} (\mathbb{R}^{1+1} )} & \lesssim \langle k\rangle^{1/4} \|\Box_k u_0\|_{L^1(\mathbb{R} )},  \label{fre-local-1d5}\\
\|\Box_k \mathscr{A} f\|_{ L^4_{x,t} (\mathbb{R}^{1+1} )} & \lesssim \langle k\rangle^{1/4} \|\Box_k f\|_{L^1_{x,t} (\mathbb{R}^{1+1})}. \label{fre-local-1d6}
\end{align}
So, using the same way as in the proof of Theorem \ref{thg1}, one can prove the result of Theorem \ref{thm1dm=3} and the details are omitted.

\section{On hyperbolic Schr\"odinger map}

In this section we prove our Corollary \ref{thmmap}.

\begin{lem} \label{Schmapa}
Let $\kappa \geq 0$, $s_1, s_2 \in M^\kappa_{r,1}$ with $r\in [1,\infty]$ and $s=(s_1,s_2,s_3)\in \mathbb{S}^2$. Suppose that $\|s_i\|_{M^\kappa_{r,1}} \leq \eta \ll 1$ for $i=1,2$. Then we have $|s_3|-1 \in M^\kappa_{r,1}$ and $\||s_3| - 1\|_{M^\kappa_{r,1}} \leq C_0 \eta$.
\end{lem}
{\bf Proof.}  It is known that $M^\kappa_{r,1} \subset L^\infty$ is a Banach algebra, we see that $\|s_i\|_{L^\infty} \ll 1$ for $i=1,2$. In view of Taylor's expansion we have
\begin{align}
|s_3|-1 & =
  \sqrt{1-s_1^2-s^2_2} -1  \nonumber\\
  & = \sum^\infty_{n=1} \frac{(-1)^n}{n \,!} (s^2_1\partial_x + s^2_2 \partial_y)^n \sqrt{1+x+y}\, \big |_{(0,0)} \nonumber\\
& = \sum^\infty_{n=1} \frac{(-1)^n}{n \,!} \sum^n_{j=0} C^j_n s^{2j}_1 \, s^{2(n-j)}_2 \prod^{n-1}_{\lambda=0}(\frac{1}{2}-\lambda).
\label{schmap1}
\end{align}
In view of the algebra property of $M^\kappa_{r,1}$,
\begin{align}
\||s_3|-1\|_{M^\kappa_{r,1}}
  &  \leq  \sum^\infty_{n=1}   \sum^n_{j=0} C^j_n \| s^{2j}_1 \, s^{2(n-j)}_2 \|_{M^\kappa_{r,1}} \nonumber\\
  &  \leq  \sum^\infty_{n=1}   \sum^n_{j=0} C^j_n C^{2n} \| s_1\|^{2j}_{M^\kappa_{r,1}} \, \|s_2 \|^{2(n-j)}_{M^\kappa_{r,1}} \nonumber\\
 &  \leq  \sum^\infty_{n=1}   \sum^n_{j=0} C^j_n C^{2n}\eta^{2n} := C_0 \eta, \label{schmap2}
\end{align}
the result follows. $\hfill\Box$

\begin{lem} \label{Schmapb}
Let $\kappa \geq 0$, $s_1, s_2 \in M^\kappa_{r,1}$ with $r\in [1,\infty]$ and $s=(s_1,s_2,s_3)\in \mathbb{S}^2$. Suppose that $\|s_i\|_{M^\kappa_{r,1}} \leq \eta \ll 1$ for $i=1,2$. Then we have $u_0:=(s_1+is_2)/(1+s_3)  \in M^\kappa_{r,1}$ and $\|u_0\|_{M^\kappa_{r,1}} \leq C_1 \eta$.
\end{lem}
{\bf Proof.} We may assume that $s_3\geq 0$. Taking $\tilde{s}_3=s_3-1$, we have $u_0:=(s_1+is_2)/2(1+\tilde{s}_3/2)$.  Let us observe that
\begin{align}
\frac{s_1+is_2}{2(1+\tilde{s}_3/2)}
  & =(s_1+is_2) \sum^\infty_{j=0}  (-1)^j  \left(\frac{\tilde{s}_3}{2} \right)^j.
\label{schmap3}
\end{align}
Using the algebra property on $M^\kappa_{r,1}$, analogous to Lemma \ref{Schmapa}, we can obtain the result, as desired. $\hfill \Box$

\medskip

By Lemmas \ref{Schmapa} and \ref{Schmapb}, we see that $u_0 =(s_1(0)+i s_2(0)) /(1+s_3(0)) \in M^2_{1,1}$ is small enough if $s_1(0), s_2(0)\in M^2_{1,1}$ with $s_0 = (s_1(0), s_2(0), s_3(0)) \in \mathbb{S}^2$ are sufficiently small. Hence, in view of Theorem \ref{thg1}, we obtain that \eqref{Schmap} has a unique solution $u\in C(\mathbb{R}, M^2_{2,1}) \cap C(\mathbb{R}, M^{3/2}_{1,1}) \cap X$.  Taking
$$
s=\left(\frac{2Re \ u}{1+|u|^2}, \  \frac{2Im \ u}{1+|u|^2}, \ \frac{|u|^2-1}{1+|u|^2}\right)
$$
and applying the same way as in Lemmas \ref{Schmapa} and \ref{Schmapb}, we have
$$
s_1, s_2, |s_3|-1 \in  C(\mathbb{R}, M^2_{2,1}) \cap C(\mathbb{R}, M^{3/2}_{1,1}).
$$
Finally, we show that $s_1, s_2, |s_3|-1 \in X$ and we need the following

\begin{lem} \label{Schmapc}
We have for $\bar{x}=(x_2,...,x_n)$,
\begin{align}
& \sum_{k\in \mathbb{Z}^n} \langle k\rangle^s \|\Box_k(u_1...u_N)\|_{L^p_{x_1}L^{q}_{\bar{x},t}} \nonumber\\
& \leq C^N \sum^N_{i=1}\left( \sum_{k^{(i)}\in \mathbb{Z}^n} \langle k^{(i)}\rangle^s \|\Box_k u_i\|_{L^p_{x_1}L^{q}_{\bar{x},t}}\right) \prod_{j\neq i, \, 1\le j\le N} \left(\sum_{k^{(j)}\in \mathbb{Z}^n} \langle k^{(j)}\rangle^s \|\Box_{k^{(j)}} u_j\|_{L^\infty_{x,t}}\right).
\label{schmap4}
\end{align}
\end{lem}
{\bf Proof.} The result was essentially obtained in \cite{WaHe07}. $\hfill\Box$

\medskip

By Taylor's expansion
$$
s_1= \sum^\infty_{j=0} (-1)^j |u|^{2j} 2Re \, u,
$$
$$
\partial_{x_1} s_1= 2Re \, u_{x_1}+ \sum^\infty_{j=1} (-1)^j |u|^{2j} 2Re \, u_{x_1} + \sum^\infty_{j=1} (-1)^j j |u|^{2j-2} 2Re \, u \partial_{x_1} |u|^2.
$$

By Lemmas \ref{Schmapc} and \ref{boxkest}, we have
$$
\sum_{\alpha=0,1}\sum_{i=1,2}\|\partial^\alpha_{x_i} s_1 \|^{\rm max \, \cap  \, str2 \, \cap\,  gstr \cap \, ant}\lesssim \delta.
$$

Finally, it suffices to estimate  $\|\partial^\alpha_{x_i} s_1 \|^{\rm sm2}$, say, we bound $\|\partial_{x_1} s_1 \|^{\rm sm2}_1$.   We have
\begin{align}
\|\partial_{x_1} s_1 \|^{\rm sm2}_1 \leq 4 \| u_{x_1}\|^{\rm sm2}_1 + \sum^\infty_{j=1} \| |u|^{2j} 2Re \, u_{x_1}\|^{\rm sm2}_1 + \sum^\infty_{j=1} j \||u|^{2j-2} 2Re \, u \partial_{x_1} |u|^2\|^{\rm sm2}_1.
\end{align}
We estimate
\begin{align}
   \| |u|^{2j} 2Re \, u_{x_1}\|^{\rm sm2}_1 & \lesssim \sum_{k\in \mathbb{Z}^2, \, |k_1|\ge |k_2|\vee 20} \left(\sum_{k^{(1)},...,k^{(j-2)}\in \mathbb{A}^{j-2, 1}_{\rm lo} }+ \sum_{k^{(1)},...,k^{(j)}\in \mathbb{A}^{j-2, 1}_{\rm hi} } \right) \langle k_1\rangle^{3/2} \nonumber \\
  & \quad\quad \times \left\|\Box_k \left(\prod^j_{s=1}\Box_{k^{(s)}} u  \prod^{2j}_{s=j+1}\Box_{k^{(s)}} \bar{u} \Box_{k^{(j-2)}}Re \, u_{x_1}  \right) \right\|_{L^\infty_{x_1}L^2_{x_2,t}} \nonumber\\
   & := \Lambda_{\rm lo} + \Lambda_{\rm hi}.
\end{align}
Using the fact
$$
\|\Box_k f\|_{L^\infty_{x_1}L^2_{x_2,t}} \le \|\Box_k f\|_{L^2_{x_2,t}L^\infty_{x_1}} \le \|\Box_k\|_{L^2_{x,t}},
$$
and in view of Lemma \ref{boxkest} and H\"older's inequality, we have
\begin{align}
  \Lambda_{\rm lo}   & \lesssim  C^j \sum_{k^{(1)},...,k^{(j-2)}\in \mathbb{Z}^{2}}
   \prod^{2j}_{s=3} \|\Box_{k^{(s)}} u\|_{L^\infty_{x,t}}  \| \Box_{k^{(j-2)}} u_{x_1} \|_{L^4_{x,t}}  \| \Box_{k^{(1)}}  u  \|_{L^4_{x,t}} \| \Box_{k^{(2)}}  u  \|_{L^4_{x,t}} \nonumber\\
    & \lesssim  C^j \|u_{x_1}  \|^{\rm str2}  (\|u\|^{\rm str2})^{2j}.
    \end{align}
Again,  Lemma \ref{boxkest} and H\"older's inequality yield
\begin{align}
  \Lambda_{\rm hi}
    & \lesssim  C^j (\|u_{x_1}  \|^{\rm sm2}  (\|u\|^{\rm str2})^{2j}+ \|u   \|^{\rm sm2} \|u_{x_1}  \|^{\rm str2} (\|u\|^{\rm str2})^{2j-1}).
    \end{align}
Collecting the estimates as in the above,  we obtain that
\begin{align}
   \| |u|^{2j} 2Re \, u_{x_1}\|^{\rm sm2}_1   \lesssim  \sum^\infty_{j=1} C^j \delta^{j-2} \lesssim \delta^3.
\end{align}
So, we have shown   $\|\partial^\alpha_{x_i} s_1 \|^{\rm sm2} \lesssim \delta$. Similarly, we have the desired estimates for $s_2$  and $|s_3|-1$. This finishes the proof of Corollary \ref{thmmap}.

\section{Initial data in weighted Sobolev spaces}

If we can show that $H^{s+b, b}(\mathbb{R}^2) \subset M^s_{1,1} (\mathbb{R}^2)$ for any $b>1$, then we get an exact proof of Corollaries  \ref{thmws1}.

\begin{prop}
Let $s\in \mathbb{R}$, $b>n/2$. We have $H^{s+b, b}(\mathbb{R}^n) \subset M^s_{1,1} (\mathbb{R}^n)$.
\end{prop}
{\bf Proof.} It suffices to consider the case $s=0$.  We have
$$
\|f\|_{M^0_{1,1}} \lesssim \left(\sum_{k\in \mathbb{Z}^n} \langle k\rangle^{-2b}\right)^{1/2} \|f\|_{M^b_{1,2}} \lesssim \|f\|_{M^b_{1,2}}.
$$
For any $\tilde{b}>n/2$, using similar way as in Lemma \ref{connorm},
\begin{align}
 \|\Box_k f\|_1
 & \lesssim   \sum_{|l|_\infty\le 1} \|\mathscr{F}^{-1}(\sigma_{k+l} \langle \xi \rangle^{-b})\|_1   \|\mathscr{F}^{-1} \sigma_{k} \langle k\rangle^{b} \mathscr{F} f\|_1 \nonumber\\
  & \lesssim    \langle k\rangle^{-b}    \|  \sigma_{k} \langle \xi\rangle^{b} \widehat{f}\|_{H^{\tilde{b}}}.
  \label{incl}
\end{align}
It  follows that
\begin{align}
\|f\|_{M^b_{1,2}} \le \left( \sum_{k\in \mathbb{Z}^n} \|  \sigma_{k} \langle \xi\rangle^{b} \widehat{f}\|^2_{H^{\tilde{b}}} \right)^{1/2}.  \label{gtosi0}
\end{align}
If $\tilde{b} \in \mathbb{N}$,  we see that
$$
\|  \sigma_{k}   \widehat{g}\|_{H^{\tilde{b}}} \lesssim \sum_{|\alpha|\le \tilde{b}} \|\partial^\alpha \widehat{g}\|_{L^2([k-1,k+1]^n)}.
$$
Hence,
\begin{align}
 \left( \sum_{k\in \mathbb{Z}^n} \|  \sigma_{k}  \widehat{g}\|^2_{H^{\tilde{b}}} \right)^{1/2} \lesssim \sum_{|\alpha|\le \tilde{b}} \|\partial^\alpha \widehat{g}\|_{2} \sim \|\langle x\rangle^{\tilde{b}} g\|_2. \label{gtosi}
\end{align}
Considering the map $T: g\to \{\sigma_k \widehat{g}: \ k\in \mathbb{Z}^n\}$,  \eqref{gtosi} implies that $T: L^2(\mathbb{R}^n:   \langle x\rangle^{\tilde{b}}dx) \to \ell^2(H^{\tilde{b}}(\mathbb{R}^n))$ is bounded for any $\tilde{b} \in \mathbb{N}\cup \{0\}$.  For any $b>n/2$, we can choose $\tilde{b}>b$ with $\tilde{b}\in \mathbb{N}$. In view of the real interpolation theory, we can interpolate $\ell^2(H^{b} (\mathbb{R}^n))$ between  $\ell^2(L^2 (\mathbb{R}^n))$ and $\ell^2(H^{\tilde{b}}(\mathbb{R}^n))$ and show that
\begin{align}
 \left( \sum_{k\in \mathbb{Z}^n} \|  \sigma_{k}  \widehat{g}\|^2_{H^{ b }} \right)^{1/2} \lesssim \sum_{|\alpha|\le \tilde{b}} \|\partial^\alpha \widehat{g}\|_{2} \sim \|\langle x\rangle^{ b } g\|_2. \label{gtosi1}
\end{align}
Taking $g= \mathscr{F}^{-1} \langle \xi\rangle^{ b }\mathscr{F}f $, we have
from \eqref{gtosi0} and \eqref{gtosi1} that
$$
\|f\|_{M^b_{1,2}} \lesssim  \|\langle x\rangle^{ b } \mathscr{F}^{-1} \langle \xi\rangle^{ b }\mathscr{F}f\|_2.
$$
Hence, we have  $\|f\|_{M^0_{1,1}} \lesssim  \| f\|_{H^{b,b}}.$  $\hfill\Box$

\section{Ill-posedness}

In this section we apply the idea as in \cite{Bou97} to show that
\begin{align}
{\rm i} u_t -  \Delta_\pm u = \overrightarrow{\lambda} \cdot\nabla(|u|^{2\kappa} u), \quad u(0,x)=u_0(x)
\label{nonls}
\end{align}
is ill-posed in $M^s_{2,1}$ if $s<1/2\kappa$. We can assume that the first coordinate of $\overrightarrow{\lambda}$ is not $0$.  Let $\varphi: \mathbb{R}^n \to [0,1]$ be a smooth function with ${\rm supp} \ \varphi \subset \{\xi\in \mathbb{R}^n:  |\xi|\le 1\}$ and  $ \varphi (\xi)=1 $ for $ \xi \in  \{\xi\in \mathbb{R}^n:  |\xi|\le 1/2\}$. Put for $0<\varepsilon \ll 1$,
\begin{align}
\widehat{u}_{0,N} & = \frac{1}{N^s}\left(\varphi\left( \varepsilon^{-1}(\xi_1-N) \right) + \varphi\left( \varepsilon^{-1}(\xi_1+N) \right)\right) \varphi\left( \varepsilon^{-1} \xi_2 \right)... \varphi\left(\varepsilon^{-1} \xi_n \right) \nonumber\\
&: = \frac{1}{N^s}(h^+_{N} + h^-_{N}). \nonumber
\end{align}
In order to show that the solution map $u_0\to u$ is not $C^{2\kappa+1}$, it suffices to prove that
$$
\sup_{t\in [0,T]} \left\|\mathscr{A} (\overrightarrow{\lambda}_1\cdot\nabla(|u|^{2\kappa} u)
 ) \right\|_{M^s_{2,1}} \lesssim \|u\|^{2\kappa+1}_{M^s_{2,1}}
$$
does not hold for $v=S(t) u_{0,N}$ if $N\gg 1$.  It is easy to see that
$$
\|u_{0,N}\|_{M^s_{2,1}} \sim_{\varepsilon} 1.
$$
Let us write
$$
\mathscr{A} (\overrightarrow{\lambda}_1\cdot\nabla(|v|^{2\kappa} v)): = c \mathscr{A} \partial_{x_1} (|v|^{2\kappa} v) + R(t,x):= M(t,x)+ R(t,x).
$$
From the argument below we will see that $M(t,x)$ contributes the main part. Denote $\xi= (\xi^{(1)},..., \xi^{(2\kappa+1)})$ for $ \xi^{(j)} \in \mathbb{R}^{n}$.  We have
\begin{align}
\widehat{M} (t, \xi^{(2\kappa+1)}) & =   e^{{\rm i} t |\xi^{(2\kappa+1)}|^2_\pm} \xi^{(2\kappa+1)}_1 \int_{\mathbb{R}^{2\kappa n}} \frac{e^{{\rm i}t P(\xi)} -1}{P(\xi)} \widehat{u}_{0,N}(\xi^{(1)})... \widehat{u}_{0,N}(\xi^{(2\kappa)}) \times \nonumber\\
& \ \ \ \  \widehat{u}_{0,N}(\xi^{(2\kappa+1)}-\xi^{(1)}-...-\xi^{(2\kappa)}) d \xi,
\label{nonlsaa}
\end{align}
where
\begin{align}
P(\xi ) = &  -\sum^{\kappa}_{j=1} |\xi^{(j)}|^2_\pm -  |\xi^{(2\kappa+1)}|^2_\pm
   +  \sum^{2\kappa}_{j=\kappa+1} |\xi^{(j)}|^2_\pm  + |\xi^{(2\kappa+1)}- \xi^{(1)}-...- \xi^{(2\kappa)}|^2_\pm  \nonumber \\
 := &  P(\bar{\xi} ) -\sum^{\kappa}_{j=1} |\xi^{(j)}_1|^2  -  |\xi^{(2\kappa+1)}_1|^2
    + \sum^{2\kappa}_{j=\kappa+1} |\xi^{(j)}_1|^2   + |\xi^{(2\kappa+1)}_1- \xi^{(1)}_1-...- \xi^{(2\kappa)}_1|^2   \nonumber
\end{align}
and $\bar{\xi}= (\bar{\xi}^{(1)},..., \bar{\xi}^{(2\kappa+1)})$ for $\bar{\xi}^{(j)}=(\xi^{(j)}_2,..., \xi^{(j)}_n)$.  We can further write
\begin{align}
\widehat{M} (t, \xi^{(2\kappa+1)}) & = \frac{1}{N^{s(2\kappa+1)}}  e^{{\rm i} t |\xi^{(2\kappa+1)}|^2_\pm} \xi^{(2\kappa+1)}_1 \int_{\mathbb{R}^{2\kappa n}} \frac{e^{{\rm i}t P(\xi )} -1}{P(\xi )} h^+_{N}(\xi^{(1)})... h^+_{N}(\xi^{(\kappa)}) \times \nonumber\\
& \ \ \ \  h^-_{N}(\xi^{(\kappa+1)})... h^-_{N}(\xi^{(2\kappa)}) h^-_{N}(\xi^{(2\kappa+1)}-\xi^{(1)}-...-\xi^{(2\kappa)}) d \xi  \nonumber\\
& \ \ \ \ + \frac{1}{N^{s(2\kappa+1)}}  e^{{\rm i} t |\xi^{(2\kappa+1)}|^2_\pm} \xi^{(2\kappa+1)}_1 \int_{\mathbb{R}^{2\kappa n}} \frac{e^{{\rm i}t P(\xi )} -1}{P(\xi )} h^-_{N}(\xi^{(1)})... h^-_{N}(\xi^{(\kappa)}) \times \nonumber\\
& \ \ \ \  h^+_{N}(\xi^{(\kappa+1)})... h^+_{N}(\xi^{(2\kappa)}) h^+_{N}(\xi^{(2\kappa+1)}-\xi^{(1)}-...-\xi^{(2\kappa)}) d \xi  +  \widehat{M}_R (t, \xi^{(2\kappa+1)}) \nonumber\\
& \ \ \ \ := \widehat{M}_1 (t, \xi^{(2\kappa+1)})+ \widehat{M}_2 (t, \xi^{(2\kappa+1)}) + \widehat{M}_R (t,  \xi^{(2\kappa+1)}).
\label{nonls2}
\end{align}
For convenience, we write for $\lambda=(\lambda_1,...,\lambda_{2\kappa+1})$, $\lambda_j=+, -$,
$$
A_\lambda  = {\rm supp} \  h^{\lambda_1}_{N}(\xi^{(1)})... h^{\lambda_{2\kappa}}_{N}(\xi^{(2\kappa)}) h^{\lambda_{2\kappa+1}}_{N}(\xi^{(2\kappa+1)}-\xi^{(1)}-...-\xi^{(2\kappa)}).
$$
We estimate $\widehat{M}_1 (t, \xi^{(2\kappa+1)})$.  By changing variables $\xi^{(j)}_1= \eta_j + N$ for $j=1,...,\kappa$ and  $\xi^{(j)}_1= \eta_j - N$ for $j=\kappa+1,..., 2\kappa+1$, one sees that
$$
P(\xi ) = O(\varepsilon)  \ \ {\rm in} \ \  A_\lambda
$$
for $\lambda_1=...=\lambda_\kappa=+$, $\lambda_{\kappa+1}=...= \lambda_{2\kappa+1}=-$. Hence, for some $\varepsilon_1>0$,
$$
|\widehat{M}_1 (t, \xi^{(2\kappa+1)}) \chi_{[-N-\varepsilon_1, -N+\varepsilon_1]\times [-\varepsilon_1,  \varepsilon_1]^{n-1}} |  \gtrsim   N^{1-s(2\kappa+1)}, \ \ t\in [T/2, T].
$$
The estimate of $\widehat{M}_2 (t, \xi^{(2\kappa+1)})$ is analogous to $\widehat{M}_1 (t, \xi^{(2\kappa+1)})$. So, we have
$$
\|M_1(t)+M_2(t)\|_{M^s_{2,1}} \gtrsim N^{1-2\kappa s}, \ \ N \gg 1, \ \ t \in [T/2, T].
$$
Let us assume that $\widehat{M}_R$ is the summation of $\widehat{M}_{R\lambda}$. In the following we show that either
$$
\|M_{R\lambda}\|_{M^s_{2,1}} \ll  N^{1-2\kappa s},
$$
or the support set of  $\widehat{M}_{R\lambda}$ is disjoint with ${\rm supp}   \widehat{M}_1 \cup {\rm supp} \widehat{M}_2 $.  We divide the estimate of $\widehat{M}_{R\lambda} (t,  \xi^{(2\kappa+1)})$ into the following three cases.

{\it Case 1.} We consider the case $\xi^{(1)}_1+...+ \xi^{(2\kappa)}_1=O(\varepsilon)$ in $A_\lambda$. For example, we estimate
\begin{align}
\widehat{M}_{R1} (t, \xi^{(2\kappa+1)}) & = \frac{ \xi^{(2\kappa+1)}_1}{N^{s(2\kappa+1)}}  e^{{\rm i} t |\xi^{(2\kappa+1)}|^2_\pm} \int_{\mathbb{R}^{2\kappa n}} \frac{e^{{\rm i}t P(\xi)} -1}{P(\xi)} h^-_{N}(\xi^{(1)}) h^+_{N}(\xi^{(2)})... h^+_{N}(\xi^{(\kappa+1)}) \times \nonumber\\
& \ \ \ \    h^-_{N}(\xi^{(\kappa+2)})... h^-_{N}(\xi^{(2\kappa)}) h^-_{N}(\xi^{(2\kappa+1)}-\xi^{(1)}-...-\xi^{(2\kappa)}) d \xi.
\label{nonls3}
\end{align}
Making a change of variables $\xi^{(j)}_1= \eta_j + N$ for $j=2,...,\kappa+1$ and  $\xi^{(j)}_1= \eta_j - N$ for $j=1, \kappa+2,..., 2\kappa+1$, we see that
$$
P(\xi ) = 2N(\eta_1+\eta_{\kappa+1}) + O(\varepsilon) \ \ {\rm in} \ \ A_\lambda.
$$
By considering  $|\eta_1+\eta_{\kappa+1}| \le 1/\sqrt{N}$ and  $|\eta_1+\eta_{\kappa+1}| \ge 1/\sqrt{N}$ in $\mathbb{R}^{2\kappa n}$, we see that
$$
|\widehat{M}_{R1} (t, \xi^{(2\kappa+1)})   |  \lesssim    N^{1/2-s(2\kappa+1)}, \ \ t\in [0, T].
$$
Noticing that ${\rm supp}\widehat{M}_{R1} \subset [-N-1/2, -N+1/2]\times [-1/2,1/2]^{n-1}$ for $0<\varepsilon \ll 1$, we have
$$
\| M_{R1} \|_{M^s_{2,1}}  \lesssim    N^{1/2- 2\kappa s }, \ \ t\in [0, T].
$$
{\it Case 2.} We consider the case $\xi^{(1)}_1+...+ \xi^{(2\kappa)}_1= \pm 2N+ O(\varepsilon)$ in $A_\lambda$. Say, we estimate
\begin{align}
\widehat{M}_{R2} (t, \xi^{(2\kappa+1)}) & = \frac{1}{N^{s(2\kappa+1)}}  e^{{\rm i} t |\xi^{(2\kappa+1)}|^2_\pm} \xi^{(2\kappa+1)}_1 \int_{\mathbb{R}^{2\kappa n}} \frac{e^{{\rm i}t P(\xi)} -1}{P(\xi)} h^+_{N}(\xi^{(1)})... h^+_{N}(\xi^{(\kappa+1)}) \times \nonumber\\
& \ \ \ \    h^-_{N}(\xi^{(\kappa+2)})... h^-_{N}(\xi^{(2\kappa)}) h^-_{N}(\xi^{(2\kappa+1)}-\xi^{(1)}-...-\xi^{(2\kappa)}) d \xi.
\label{nonls4}
\end{align}
Changing variables $\xi^{(j)}_1= \eta_j + N$ for $j=1,...,\kappa+1, 2\kappa+1$ and  $\xi^{(j)}_1= \eta_j - N$ for $j=\kappa+2,..., 2\kappa$, we have
$$
P(\xi ) = 4N(\eta_{\kappa+1} + \eta_{2\kappa+1}) + O(\varepsilon)  \ \ {\rm in} \ \  A_\lambda
$$
Similarly as in Case 1, we have
$$
\| M_{R2} \|_{M^s_{2,1}}  \lesssim    N^{1/2- 2\kappa s }, \ \ t\in [0, T].
$$
{\it Case 3.} We consider the case $\xi^{(1)}_1+...+ \xi^{(2\kappa)}_1= \pm 2kN+ O(\varepsilon)$ in $A_\lambda$ with $k\ge 2$. We easily see that the support sets of  $\widehat{M}_{R\lambda}$ in this case never overlap with the support sets of $\widehat{M}_1$  and $\widehat{M}_2$.\\

Consider the following problem
\begin{align}
{\rm i} u_t -  \Delta_\pm u =  |\partial_{x_1} u|^{2\kappa} \partial_{x_1} u , \quad u(0,x)=u_0(x).
\label{nonlsa}
\end{align}
Taking $v= \partial_{x_1} u$, we see that $v$ satisfies
\begin{align}
{\rm i} v_t -  \Delta_\pm v = \partial_{x_1} (|v|^{2\kappa} v) , \quad  v(0,x)=v_0(x).
\label{nonlsab}
\end{align}
Using the same way as in the above, one sees that \eqref{nonlsa} is ill-posed in $M^{s}_{2,1}$ if $s<1+1/2\kappa$.

\medskip

\noindent {\bf Remark.} From the proof above we see that
\begin{align}
{\rm i} v_t -  \Delta  v = \mu |v|^{2\nu} v  , \quad  v(0,x)=v_0(x)
\label{nonlsabcns}
\end{align}
is ill posed in $M^s_{2,1}$ if $s<0$, i.e., the solution map is not $C^{2\nu +1}$ from $M^s_{2,1}$ into $C([0,T]; M^s_{2,1})$ for any $s<0$ and $T>0$. This implies that the well posed results in $M^0_{2,1}$ for NLS \eqref{nonlsabcns}  obtained in \cite{WaZhGu06, WaHe07} are also sharp with respect the spatial regularity index.

\begin{appendix}

\section{Gabor frame}

We collect some results used in this paper for the Gabor frame, see for instance, Gr\"ochenig \cite{Gr01}. Gabor frame is a fundamental tool in the theory of time-frequency analysis, which was first proposed by Gabor \cite{Ga46} in  1946.  A system $\{e_j: \ j\in J\}$ in a Hilbert space $\mathcal{H}$ is said to be a frame if there exists two positive constant $A,B>0$ such that for all $f\in \mathcal{H}$,
$$
A\|f\| \le \left(\sum_{j\in J} |\langle f, e_j\rangle|^2 \right)^{1/2} \le B\|f\|.
$$
For convenience, we write $T_x f (y) = f(y-x)$, $M_\xi f (y)= e^{{\rm i} y\cdot \xi} f(y)$, $T_x$ is a translation by $x$ and $M_\xi$ is a modulation by $\xi$.  Let $g\in L^2(\mathbb{R}^n)$ and $\alpha,\beta>0$. If
$$
\mathcal{G}(g,\alpha,\beta):= \{T_{\alpha l} M_{\beta k} g: \ \ k, l\in \mathbb{Z}^n\}
$$
is a frame in $L^2$, then it is said to be a Gabor frame in $L^2$.

\begin{prop} \label{gabor1}
Let $\mathcal{G}(g,\alpha,\beta)$ be a Gabor frame in $L^2$. Then any $f\in L^2$ has an expansion
$$
f= \sum_{k,l\in \mathbb{Z}^n} c_{kl} T_{\alpha l} M_{\beta k} g.
$$
Moreover, $\|f\|_2 \sim \|(c_{kl})\|_{\ell^2}$.
\end{prop}
Unfortunately, the generalization of Gabor frame in $L^p$ with $p\neq 2$ is not available and the Gabor expansion only holds for the case $p=2$. However, Proposition \ref{gabor1} also holds for modulation spaces $M^s_{p,q}$ with $1\le p,q<\infty$:
\begin{prop} \label{gabor2}
Let $\mathcal{G}(g,\alpha,\beta)$ be a Gabor frame in $L^2$. Then any $f\in M^s_{p,q}$ has an expansion
$$
f= \sum_{k,l\in \mathbb{Z}^n} c_{kl} T_{\alpha l} M_{\beta k} g.
$$
Moreover, $\|f\|_{M^s_{p,q}} \sim \|(c_{kl})\|_{m^s_{p,q}}$, where
$$
\|(c_{kl})\|_{m^s_{p,q}} = \left(\sum_{k}\langle k\rangle^{sq} \left(\sum_{l  }| c_{kl} |^p\right)^{q/p} \right)^{1/q}.
$$
\end{prop}
A basic example of the Gabor frame is $\{e^{{\rm i} kx}e^{-|x-l|^2/2}: \ k,l\in \mathbb{Z}^n\}$.

\section{Function-sequence convolution}

 Considering the Gabor frame expression for the solutions of linear Schr\"odinger equations, we need to treat the convolution on variables $x\in \mathbb{R}^n$ and $l \in \mathbb{Z}^n$.  Since $x$ and $l$ belong to different measure spaces, we can not directly use Young's and Hardy-Littlewood-Sobolev's inequalities.  So, we need the following

\begin{lem}\label{appb}
Let $1\leq p, r \leq \infty$. Assume that $\theta>0$, $\theta > 1/r'+ 1/p$ with $p\ge r$, or $\theta = 1/r'+1/p \in (0,1)$ and $1<p<\infty$.
Then we have for any $b,c\in \mathbb{R}$ with $|c|\geq 1$,
\begin{align} \label{convol}
  \left\| \sum_{l \in \mathbb{Z} }  |a_{l}|    \left(1+  \frac{| x - l+ b|}{ |c|}\right)^{-\theta} \right\|_{L^{ p }_{ x }} \lesssim \langle c\rangle^{1/p+1/r'} \|(a_l)\|_{\ell^r}.
\end{align}
\end{lem}
\noindent{\bf Proof.} For convenience, we denote by $[x]$ the integer part of $x\in \mathbb{R}$. We have
\begin{align}
  & \sum_{l \in \mathbb{Z} }  |a_{l}|    \left(1+  \frac{| x - l+ b|}{ |c|}\right)^{-\theta} \nonumber\\
  &=   \sum_{l \in \mathbb{Z} } \int^{l+1}_{l}  |a_{l}|    \left(1+  \frac{| x - l+ b|}{ |c|}\right)^{-\theta} dy \nonumber\\
    &=   \left(\sum_{l \geq [x+b]+2} + \sum_{l \leq [x+b]-1}\right) \int^{l+1}_{l}  |a_{l}|    \left(1+  \frac{| x - l+ b|}{ |c|}\right)^{-\theta} dy \nonumber\\
    & \quad + \sum^{[x+b]+1}_{l=[x+b]}\int^{l+1}_{l} |a_{l}|    \left(1+  \frac{| x  + b-l|}{ |c|}\right)^{-\theta} dy.
  \end{align}
Noticing that $|x-l+b |\geq |x-y+b+1|$ for any $y\in [l, l+1)$, $l\ge [x+b]+2$,  we have
\begin{align}
      & \sum_{l \geq [x+b]+2}   \int^{l+1}_{l}  |a_{l}|    \left(1+  \frac{| x - l+ b|}{ |c|}\right)^{-\theta} dy \nonumber\\
    &  \leq  \sum_{l \geq [x+b]+2} \int^{l+1}_{l} |a_{l}|    \left(1+  \frac{| x  -y+b +1|}{ |c|}\right)^{-\theta} dy \nonumber\\
    &  \leq  \int_{\mathbb{R}}    \left(1+  \frac{| x  -y+b +1|}{ |c|}\right)^{-\theta} \sum_{l \in \mathbb{Z}} |a_{l}| \chi_{[l, l+1)} (y) dy.
  \end{align}
Similarly,
\begin{align}
      & \sum_{l \leq [x+b]-1}   \int^{l+1}_{l}     \left(1+  \frac{| x - l+ b|}{ |c|}\right)^{-\theta} |a_{l}|  dy \nonumber\\
    &  \leq  \int_{\mathbb{R}}    \left(1+  \frac{| x  -y+b|}{ |c|}\right)^{-\theta} \sum_{l \in \mathbb{Z}} |a_{l}| \chi_{[l, l+1)} (y) dy,
  \end{align}
and noticing that $|c| \geq 1$, for $l=[x+b], [x+b]+1$ one has that
\begin{align}
      &     \int^{l+1}_{l}     \left(1+  \frac{| x + b-l|}{ |c|}\right)^{-\theta} |a_{l}|  dy \nonumber\\
      &  \leq  4^\theta \int_{\mathbb{R}}    \left(1+  \frac{| x  -y+b|}{ |c|}\right)^{-\theta}   |a_{l}| \chi_{[l, l+1)} (y) dy \nonumber\\
    &  \leq  4^\theta \int_{\mathbb{R}}    \left(1+  \frac{| x  -y+b|}{ |c|}\right)^{-\theta} \sum_{l \in \mathbb{Z}} |a_{l}| \chi_{[l, l+1)} (y) dy.
  \end{align}
Applying Young's inequality in the case $\theta >1/r'+1/p$, we have
\begin{align}
 &      \left\| \int_{\mathbb{R}}    \left(1+  \frac{| x  -y+b|}{ |c|}\right)^{-\theta} \sum_{l \in \mathbb{Z}} |a_{l}| \chi_{[l, l+1)} (y) dy \right\|_{L^p_x} \nonumber\\
  & \lesssim
    \|(1+|\cdot|/|c|)^{-\theta}\|_{L^{pr'/(p+r')}} \left\| \sum_{l \in \mathbb{Z}} |a_{l}| \chi_{[l, l+1)} \right\|_{L^r} \nonumber\\
   & \lesssim
    \langle c\rangle^{1/r'+1/p} \left\| \sum_{l \in \mathbb{Z}} |a_{l}| \chi_{[l, l+1)} \right\|_{L^r}.
    \end{align}
Using Hardy-Littlewood-Sobolev's inequality in the case $\theta =1/r'+1/p$, we also have the result, as desired.  $\hfill\Box$

\medskip

We note that the hidden constant in the right hand side of \eqref{convol} is independent of $b,c\in \mathbb{R}$ with $|c|\geq 1$.

\section{Blow up solution of \eqref{Schmap} in 2D}

Guo and Yang \cite{GuYa} found a class of solutions of the 2D Schr\"odinger map equation in the elliptic case, which contain a blow up solution for \eqref{Schmap} in 2D elliptic case. We can easily generalize their result to 2D hyperbolic case of \eqref{Schmap}.
Let $\langle t \rangle = (1+t^2)^{1/2}$ and
$$
s_1 = - \frac{t}{\langle t \rangle}, \ \ s_2= \frac{1}{\langle t \rangle} \sin  \frac{x^2_1-x^2_2}{4 \langle t \rangle}, \ \  s_3= \frac{1}{\langle t \rangle} \cos  \frac{x^2_1-x^2_2}{ 4 \langle t \rangle}.
$$
One easily sees that $s=(s_1, s_2, s_3)$  satisfies
$$
s_t = s\times \Box s, \ \   \Box=\partial^2_{x_1} -  \partial^2_{x_2}.
$$
Indeed, noticing that $\Box s_1=0$, it suffices to show
$$
\left(
\begin{array}{c}
\partial_t s_1\\
\partial_t s_2\\
\partial_t s_3
\end{array}
\right)
=
\left(
\begin{array}{c}
  s_2 \Box s_3 - s_3\Box s_2\\
- s_1 \Box s_3 \\
 s_1\Box s_2
\end{array}
\right).
$$
By a simple calculation, we see that
$$
\Box s_2 =  \frac{1}{\langle t \rangle^2 } \cos  \frac{(x^2_1-x^2_2)}{ 4 \langle t \rangle} - \frac{(x^2_1-x^2_2)}{ 4 \langle t \rangle^3 } \sin  \frac{(x^2_1-x^2_2)}{ 4 \langle t \rangle},
$$
$$
\Box s_3 = - \frac{1}{\langle t \rangle^2 } \sin  \frac{(x^2_1-x^2_2)}{ 4 \langle t \rangle} - \frac{(x^2_1-x^2_2)}{ 4 \langle t \rangle^3 } \cos  \frac{(x^2_1-x^2_2)}{ 4 \langle t \rangle}.
$$
Moreover, by an easy calculation to $\partial_t s_i$, we easily find that $s$ is a solution. In view of the stereographic projection, we immediately have

\begin{prop} \label{blowup}
$$
u (t,x) = \frac{-t + {\rm i}  \sin \frac{(x^2_1-x^2_2)}{ 4 \langle t \rangle}   }{\langle t \rangle  +  \cos  \frac{(x^2_1-x^2_2)}{ 4 \langle t \rangle}}
$$
is a solution of
\begin{align}
{\rm i} u_t -  \Box  u     = \frac{2\bar{u}}{1+|u|^2} (  u^2_{x_1}- u^2_{x_2}),
 \end{align}
which blows up at $t=0$.
\end{prop}
{\bf Proof.} It suffices to show that $u$ blows up at $t=0$. Indeed,
$$
u (0,x) = \frac{\sin \frac{(x^2_1-x^2_2)}{ 4 }   }{1+  \cos  \frac{(x^2_1-x^2_2)}{ 4}},
$$
one of the blow up curves is $  x^2_1-x^2_2  = 4 \pi$.  $\hfill \Box$

Replacing $t$ by $t-T $, we have a bow up solution at $t=T$:
$$
u (t,x) = \frac{T-t + {\rm i}  \sin \frac{(x^2_1-x^2_2)}{ 4 \langle t-T \rangle}   }{\langle t -T\rangle  +  \cos  \frac{(x^2_1-x^2_2)}{ 4 \langle t-T \rangle}},
$$
Since the solution has no decay as $|x|\to \infty$, it does not belong to any Sobolev or modulation spaces.

\end{appendix}

\medskip

{\bf Acknowledgment.} The author is grateful to Professors K. Gr\"ochenig and M. Sugimoto for their enlightening conversations on Gabor frame. He is also indebted to Professor Ganshan Yang for his discussions on the blow up solutions of the Schr\"odinger map equation in the elliptic case.
The author was supported in part by an NSFC grant and a 973 grant.

\end{document}